# MINIMAL REDUCTION TYPE AND AFFINE SPRINGER FIBERS

ANLONG CHUA

ABSTRACT. We extend a result of Yun on minimal reduction types to the parahoric case. This implies a uniqueness property for 2-special representations appearing in the cohomology of certain affine Springer fibers. Using this, we settle a conjecture of Lusztig on strata in a reductive group.

## 1. INTRODUCTION

1.1. **Preliminary notation.** Let $G$ be a connected, reductive group over $\mathbb{C}$ with Lie algebra $\mathfrak{g}$. Fix a Borel subgroup $B \subset G$ and maximal torus $T \subset G$. Denote their Lie algebras by $\mathfrak{b}, \mathfrak{t}$ respectively. Let $W$ be the Weyl group of $G$, and denote the set of conjugacy classes in $W$ by $\underline{W}$. We denote by $\mathrm{Irr}(W)$ the set of irreducible representations of $W$. Let $\mathcal{N} \subset \mathfrak{g}$ denote the nilpotent cone, and $\underline{\mathcal{N}}$ be the finite set of nilpotent orbits.

Let $F = \mathbb{C}((t))$ be the field of Laurent series in $t$. The subring of Taylor series $\mathcal{O} = \mathbb{C}[[t]]$ is its ring of integers with respect to the $t$-adic valuation. Let $I \subset G(F)$ be the preimage of $B$ under the evaluation map $G(\mathcal{O}) \to G(\mathbb{C})$. Let $P \subset G(F)$ be a parahoric subgroup. Denote by $\mathrm{Fl}^P = G(F)/P$ the affine partial flag variety associated to $P$. When $P = I$, we write $\mathrm{Fl}^I = \mathrm{Fl}$ and when $P = G(\mathcal{O})$ we write $\mathrm{Fl}^P = \mathrm{Gr}$.

For any $\gamma \in \mathfrak{g}(F)$, let
$$\mathrm{Fl}^P_\gamma = \{gP \in \mathrm{Fl}^P \mid \mathrm{Ad}(g^{-1})\gamma \in \mathrm{Lie}\, P\}$$
be the partial affine Springer fiber associated to $\gamma$. It comes equipped with an evaluation map
$$\mathrm{ev}_P \colon \mathrm{Fl}^P_\gamma \to [\mathfrak{g}_P/G_P]$$
sending
$$gP \mapsto \mathrm{Ad}(g^{-1})\gamma \pmod{P^+} \in \mathfrak{g}_P/G_P.$$
where $G_P$ is the Levi quotient of $P$.

1.2. **Motivation.** Recall that conjugacy classes of maximal tori in $G(F)$ are indexed by $\underline{W}$ [13, Section 1, Lemma 2]. Let $\gamma \in \mathfrak{g}(F)$ be regular semisimple. We say that it has type $[w] \in \underline{W}$ if its centralizer is of type $[w]$. Denote by $\mathfrak{g}^{\heartsuit}(F)_{[w]}$ the subset of topologically nilpotent regular semisimple elements of type $[w]$.

Let us write $\mathrm{ev}_P(\gamma)$ for the set of nilpotent orbits appearing under the evaluation map. We sometimes call these $P$-reduction types. It forms a sub-poset of the partially ordered set $\underline{\mathcal{N}}_{G_P}$.

In [27], Yun proved:

**Theorem 1.** [27, Theorem 1] *There is an open dense subset of $\mathfrak{g}^{\heartsuit}(F)_{[w]}$ such that for any $\gamma \in \mathfrak{g}(F)_{[w]}$, $\mathrm{ev}_{G(\mathcal{O})}(\gamma)$ has a unique minimal element.*

For the definition of the topology on $\mathfrak{g}^\heartsuit(F)_{[w]}$, see Section 2.3.

If one simply replaces $\mathrm{ev}_{G(\mathcal{O})}(\gamma)$ with $\mathrm{ev}_P(\gamma)$ for a different parahoric subgroup $P \subset G(F)$, the resulting statement is false (see Example 33). The goal of this paper is to explain and analyze the failure of the statement to be true for arbitrary $P$.

1.3. **Main results.** Before stating the main theorem, we require some notation. For a nilpotent orbit $O_P \in \mathrm{ev}_P(\gamma)$, let $b_{O_P}$ denote the dimension of the Springer fiber of an element in $O_P$ in the Levi quotient $G_P$. By the dimension formula for Springer fibers, a maximal value of $b_{O_P}$ corresponds to a minimal orbit $O_P$.

**Definition 2.** *Let $b^*(\mathrm{Fl}_\gamma)$ be the maximal value of $b_{O_P}$ over any $P \subset G(F)$ and any $O_P \subset \mathrm{ev}_P(\gamma)$.*

Let $W'$ be a subgroup of $W$. In this paper, $W'$ will be a so-called parahoric subgroup. Lusztig-Macdonald-Spaltenstein have defined a functor $j_{W'}^W \colon \mathrm{Irr}(W') \to \mathrm{Irr}(W)$.

**Definition 3.** *Let $\gamma \in \mathfrak{g}(F)$. We define the set*
$$R(\gamma) = \{j_{W_P}^W E_{O_P} \mid P \subset G(F), O_P \in \mathrm{ev}_P(\gamma), b_{O_P} = b^*(\mathrm{Fl}_\gamma)\}.$$

Here $E_{O_P}$ is the Springer representation (of the Levi quotient of $P$) corresponding to $O_P$.

The main result of this paper is:

**Theorem 4.** *Assume that $[w]$ is in the image of the Kazhdan-Lusztig map. There exists an open dense subset of $\mathfrak{g}(F)_{[w]}$ such that for any $\gamma \in \mathfrak{g}(F)_{[w]}$, $R(\gamma)$ has a single element.*

The statement of Theorem 4 was proposed to us by Zhiwei Yun as a conjecture. The Kazhdan-Lusztig map was defined in [13, Section 9.1]. Given a parahoric subgroup $P \subset G(F)$ with Levi quotient $P/P^+ = G_P$, and a nilpotent orbit $O_P \in \mathcal{N}_{G_P}$, Kazhdan and Lusztig define, using the topology of $\mathrm{Lie}\,P$, a conjugacy class $\overline{\mathrm{KL}_G^P}(O_P) \in \underline{W}$. See Section 2.4 for more details.

Let us spell out the consequence of Theorem 4 in more detail. Let $P \subset G(F)$ be a parahoric subgroup, and $O_P$, $O'_P$ two $P$-minimal reduction types. As remarked before, the minimal reduction type may not be unique. However, the theorem asserts that these orbits are related by the equation $j_{W_P}^W E_{O_P} = j_{W_P}^W E_{O'_P}$.

More interestingly, let $Q \subset G(F)$ be another parahoric subgroup, and $O_Q$ be a $Q$-minimal reduction type. Note that $O_P$ and $O_Q$ are orbits in different reductive groups. A priori, it is unclear what the relation between these orbits are. Theorem 4 asserts that their Springer representations satisfy $j_{W_P}^W E_{O_P} = j_{W_Q}^W E_{O_Q}$.

1.3.1. *Affine Springer representations.* We describe a translation Theorem 4 into the language of affine Springer fibers.

The set $\mathrm{Fl}_\gamma$ can be given the structure of a scheme, and $H_c^*(\mathrm{Fl}_\gamma)$ can be given the structure of a $W$-module (see Section 3.1). The centralizer $Z(\gamma) = Z_{G(F)}(\gamma)$ acts on $\mathrm{Fl}_\gamma$, and hence on $H_c^*(\mathrm{Fl}_\gamma)$. Let us denote the fixed points of this action by $H_c^*(\mathrm{Fl}_\gamma)^{st}$.

Let $C \subset \mathrm{Irr}(W)$ be a class of representations of $W$. Let us decompose into isotypic components
$$H_c^{\mathrm{top}}(\mathrm{Fl}_\gamma)^{st} = \bigoplus_{E \in \mathrm{Irr}\,W} E \otimes V_E.$$

Define
$$b_C(\text{Fl}_\gamma) = \max\{b_E \mid V_E \neq 0 \text{ and } E \in C\}$$
where $b_E$ is the $b$-invariant of $E$ defined by Lusztig in [15], see also Section 2.1.

**Definition 5.** *For a class $C \subset \text{Irr}(W)$ of representations of $W$, we say that $\text{Fl}_\gamma$ has property $P_C$ if the $W$-module $H_c^{\text{top}}(\text{Fl}_\gamma)^{st}$ has a unique isomorphism class of representation $E$ such that $b_E = b_C(\text{Fl}_\gamma)$ and $E \in C$.*

We call a representation of $W$ a Springer representation if it corresponds to a nilpotent orbit with the trivial local system under the Springer correspondence. Lusztig has defined two additional classes of irreducible representations of $W$: special and 2-special representations. We have inclusions

{2-special representations} $\supset$ {Springer representations} $\supset$ {special representations}.

Note that these classes depend on $G$ as well as $W$. For more details, see Section 2.1.

Theorem 1 can be roughly interpreted as saying that there exists an open subset $U \subset \mathfrak{g}^\heartsuit(F)_{[w]}$ such that for any $\gamma \in U$, $\text{Fl}_\gamma$ has property $P_{\text{Springer}}$.

Our main theorem implies

**Theorem 6.** *Assume that $[w]$ is in the image of the Kazhdan-Lusztig map. Then there is an open dense subset of $\mathfrak{g}^\heartsuit(F)_{[w]}$ such that for any $\gamma \in \mathfrak{g}^\heartsuit(F)_{[w]}$, $\text{Fl}_\gamma$ has property $P_{\text{2-special}}$.*

1.3.2. *Application to strata.* In [19], Lusztig defines a partition of $G$ into subsets called strata. Strata are unions of conjugacy classes of $G$, and have desirable properties. For example, unlike conjugacy classes, strata are locally closed [6, Theorem 2.2]. The definition of strata involves the Springer correspondence and the combinatorial geometry of $G/B$.

In [19, Section 6.4], Lusztig conjectures that strata can also be described in terms of (parahoric) Kazhdan-Lusztig maps, which are defined using the topology of the loop Lie algebra (see Section 2.4). Toward this conjecture, we prove (see Proposition 25):

**Theorem 7.** *Let $P$ and $Q$ be parahoric subgroups of $G(F)$. Let $O_P$ and $O_Q$ be nilpotent orbits in their respective Levi quotients. Then*
$$\text{KL}_G^P(O_P) = \text{KL}_G^Q(O_Q)$$
*if and only if*
$$j_{W_P}^W E_{O_P} = j_{W_Q}^W E_{O_Q}. \tag{1}$$

Note that (1) involves data coming from classical Springer theory, which a priori has no relation to the affine story. As a corollary of Theorem 7, we deduce:

**Corollary 8.** *Lusztig's conjecture in [19, Section 6.4] holds.*

For more details and complete statements, see Section 2.5.

1.4. **Paper overview.** We make precise the notation used in this introduction in Section 2. Next, we show the equivalence of the various forms of Theorem 4 (Theorem 6, Theorem 7) in Section 3.

In Section 4, we reduce the proof of Theorem 4 to the case $G$ is almost simple, $[w]$ is elliptic, and show that one only needs to consider one isomorphism class of $G$ in its isogeny class.

Finally, we prove the theorem in these cases in Section 5. For classical types, the proof uses *skeleta* first introduced by Yun ([27, Section 7]). In exceptional types, the proof mainly consists of checking tables, with some complications in type $E_8$.

1.5. **Acknowledgments.** I must sincerely thank Zhiwei Yun for sharing his conjecture, which resulted in this paper, as well as many inspiring discussions. I have also benefited a lot from discussions with Roman Bezrukavnikov, Ivan Losev and Yakov Varshavsky on topics appearing in this paper.

## 2. Background

In this section, we make precise the notation and objects introduced in Section 1. We also prove some propositions in preparation for the proof of Theorem 6.

2.1. **Representations of $W$.** In this subsection, we will recall the definitions of Springer, special and 2-special representations.

2.1.1. *Springer representations.* Recall that the Springer correspondence gives an injection

$$\text{Irr}(W) \hookrightarrow \{(\mathcal{O}, \mathscr{L}) \mid \mathcal{O} \in \underline{\mathcal{N}}, \mathcal{L} \text{ an irreducible } G\text{-equivariant local system on } \mathcal{O}\}.$$

The representations corresponding to the trivial local system are called *Springer representations*. In the literature, there are (at least) two constructions of this map; we use Lusztig's version [7, Section 13.3].

2.1.2. *Special representations.* Let $E$ be an irreducible representation of $W$. To $E$, Lusztig [15] attaches two integers $a_E$ and $b_E$, as follows. Let $D_E(q)$ be the generic degree polynomial associated to $E$. Then $D_E(q) = \gamma_E q^{a_E} +$ higher powers. The integer $b_E$ is defined to be the smallest nonnegative integer $i$ such that $E$ appears in the $i$th symmetric power of the reflection representation of $W$.

Lusztig observes [15, 2.1] that $a_E \leq b_E$. When this is an equality, $E$ is called *special*. It is an empirical observation that the set of special representations is contained in the set of Springer representations. Nilpotent orbits which correspond to special representations are called special.

2.1.3. *2-special representations.* Lusztig defines a larger class of representations called 2-special representations in [19]. To give the definition, we must first recall Lusztig's $j$-induction operation.

Let $W_1 \subset W$ be a subgroup, and $E_1$ be a $W_1$-module such that $E_1$ appears in the $b_E$th symmetric power of the reflection representation of $W_1$ with multiplicity 1. Let $E = \text{Ind}_{W_1}^W E_1$. Then it is known (due to Macdonald-Lusztig-Spaltenstein, [11, Theorem 5.2.6]) that for any $E'$ appearing in $E$, $b_{E'} \geq b_{E_1}$, and moreover, equality is achieved for a unique $E''$. We shall denote $j_{W_1}^W E_1 := E''$. It is known that Springer representations satisfy the condition on $E_1$ ([19, Section 1.1]).

Let $s \in G$ be a semisimple element. Then $G_s := Z_G(s)^\circ$ is a connected reductive group, sometimes called a pseudo-Levi subgroup of $G$. Let $\zeta \in G_s^\vee$ be a semisimple element in the Langlands dual of $G_s$. Then the group $(G_s^\vee)_\zeta = Z_{G_s^\vee}(\zeta)$ is sometimes called an endoscopic group for $G_s$. Let us denote its Weyl group by $W_{s,\zeta}$.

**Definition 9.** [19, Section 1.1(a)] *A 2-special representation of $G$ is one that can be written as $j_{W_{s,\zeta}}^W E_{s,\zeta}$ for $E_{s,\zeta}$ a special representation of $W_{s,\zeta}$.*

Equivalently, one may first take an endoscopic group $H$ of $G$, take a pseudolevi subgroup $H_\zeta$ of $H$, then define 2-special representations to be the ones obtained by $j$-induction from special representations of the Weyl group of $H_\zeta$.

Taking Langlands duals swaps these two equivalent definitions. It follows that:

**Lemma 10.** [19, Section 1.1(b)] *A $W$-module $E$ is 2-special for $G$ if and only if it is 2-special for $G^\vee$.*

2.2. **Parahoric subgroups.** In this paper, all parahoric subgroups are taken to be standard, i.e. containing our choice of Iwahori $I$. If $P$ is a parahoric subgroup, we denote by $P^+$ its pro-unipotent radical. There is a short exact sequence

$$1 \to P^+ \to P \to G_P \to 1.$$

The quotient $G_P$ is a reductive group over $\mathbb{C}$. We denote its Weyl group by $W_P$.

Let $\mathbb{Z}\Phi^\vee \subset \mathbb{X}_*(T)$ be the coroot lattice. Define

$$W_{\text{aff}} \coloneqq \mathbb{Z}\Phi^\vee \rtimes W, \qquad \widetilde{W} = \mathbb{X}_*(T) \rtimes W$$

to be the affine Weyl group and extended affine Weyl group of $G$. If $G$ is almost simple, $W_{\text{aff}}$ is a Coxeter group with simple reflections $\Sigma_{\text{aff}} \coloneqq \Sigma \cup \{s_0\}$ where $s_0$ is the reflection corresponding to the affine root $\alpha_0 = 1 - \beta$; here $\beta$ is the longest root. In this situation, standard parahoric subgroups are in bijection with proper subsets of $\Sigma_{\text{aff}}$.

Let $C \subset \mathbb{X}_*(T) \otimes \mathbb{R}$ be the compact affine spherical alcove. It is the region where the affine roots satisfy $\alpha_i \geq 0$ (for $0 \leq i \leq \text{rank}\, G$. Given $x \in C$, one can associate a parahoric subgroup $P_x \subset G(F)$ by declaring $P_x$ to be generated by $T(\mathcal{O})$ and $U_\alpha(\pi^n \mathcal{O})$ for all $(\alpha, n)$ such that $\alpha(x) + n \geq 0$. Every parahoric subgroup can be expressed this way. In this case, a set of simple roots for $G_P$ is $\Sigma_P = \{\alpha_i \mid \alpha_i(x) = 0, 0 \leq i \leq \text{rank}\, G\}$. Then $W_P$ is the Coxeter group generated by simple reflections in $\Sigma_{\text{aff}}$ corresponding to the roots in $\Sigma_P$. Thus, we naturally have an embedding $W_P \subset W_{\text{aff}}$. The projection

$$W_P \subset W_{\text{aff}} = \mathbb{Z}\Phi^\vee \rtimes W \to W,$$

maps $W_P$ isomorphically onto its image in $W$.

2.3. **Arc spaces.** Let $X$ be an affine scheme of finite type over $\mathbb{C}$. We denote by $LX$ the functor $LX(R) \mapsto X(R((t)))$ on the category of $\mathbb{C}$-algebras. This is representable by an ind-scheme. Similarly, define $L^+ X$ to be the functor $L^+ X(R) \mapsto X(R[[t]])$. It is representable by a scheme (not of finite type). We have the reduction morphism $\text{ev}_t \colon L^+ X \to X$ given by reduction mod $t$.

In particular, let $\text{ev}_t^{-1}(B) =: I \subset L^+ G \subset LG$. Let $I \subset P \subset LG(\mathbb{C})$ be a standard parahoric subgroup. We define the corresponding affine partial flag variety $\text{Fl}^P$ to be the sheafification of the functor $\text{Fl}^P(R) \mapsto LG(R)/P(R)$ on the category of $\mathbb{C}$-algebras. When $P = I$, we abbreviate $\text{Fl}^I = \text{Fl}$, and call this the affine flag variety of $G$.

2.3.1. *The Chevalley base.* Let $\mathfrak{c} = \mathfrak{g}//G = \mathfrak{t}//W$. Let $\chi = \mathfrak{g} \to \mathfrak{c}$ be the quotient. It induces a map $L\mathfrak{g} \to L\mathfrak{c}$, which we will also denote by $\chi$.

Let $L^\heartsuit \mathfrak{g} \subset L^+ \mathfrak{g}$ be the subscheme of topologically nilpotent, regular semisimple elements. We will use $L^\heartsuit(-)$ to denote topologically nilpotent, regular semisimple elements in various arc spaces related to $\mathfrak{g}$.

For any $[w] \in \underline{W}$, let $(L^\heartsuit \mathfrak{g})_{[w]} \subset L\mathfrak{g}$ be the subset of topologically nilpotent, regular semisimple elements of type $[w]$. Let $(L^\heartsuit \mathfrak{c})_{[w]} = \chi(L^\heartsuit \mathfrak{g})_{[w]}$ its image in $L\mathfrak{c}$. Yun shows in [27, Section 4] that $(L^\heartsuit \mathfrak{c})_{[w]}$ is a scheme with reasonable properties.

Let $\Delta \colon \mathfrak{t} \to \mathbb{A}^1$,

$$\Delta(x) = \prod_{\alpha \in \Phi} \alpha(x) \tag{2}$$

be the discriminant function. We extend this to a function on $\mathfrak{g}$ by $\mathrm{Ad}(G)$-invariance, and then to a function on $L\mathfrak{g}$ by extension of scalars.

Note that $\gamma \in L\mathfrak{g}$ is conjugate to an element in $\mathfrak{t}\left(\mathbb{C}((t^{1/m}))\right)$ after passing to the algebraic closure $\overline{\mathbb{C}((t))} = \bigcup_{k \geq 1} \mathbb{C}((t^{1/k}))$. Extending the $t$-adic valuation so that $\mathrm{val}(t^{m/n}) = \frac{m}{n}$, we can make sense of $\mathrm{val}(\Delta(\gamma))$.

Define

$$\delta(\gamma) = \frac{\mathrm{val}(\Delta(\gamma)) - (\dim \mathfrak{t} - \dim \mathfrak{t}^w)}{2} \tag{3}$$

where $\mathfrak{t}^w$ is the subspace where $w$ acts by eigenvalue 1. This is Bezrukavnikov's formula for the dimension of the affine Springer fiber $\mathrm{Fl}_\gamma$, see Section 3.1.

We now recall a definition due to Yun [27, Definition 1.6]:

**Definition 11.** *For $[w] \in W$, let $\delta_{[w]} = \min\{\delta(\gamma) \mid \gamma \in (L^\heartsuit \mathfrak{g})_{[w]}\}$. An element $\gamma \in L^\heartsuit \mathfrak{g}$ is called* shallow *of type $[w] \in \underline{W}$ if $\gamma$ is of type $[w]$, and $\delta(\gamma) = \delta_{[w]}$.*

We denote the set of shallow elements of type $[w]$ by $(L^\heartsuit \mathfrak{g})^{sh}_{[w]}$. Let $(L^\heartsuit \mathfrak{c})^{sh}_{[w]}$ be its image under $\chi$. This is a locally closed subscheme of $L^+ \mathfrak{c}$.

2.3.2. *Topological notions.* Let $X$ be an affine scheme of finite type over $\mathbb{C}$. Define, for any $n \geq 0$, $L_n^+ X$ to be the truncated arc space representing $R \mapsto X(R[[t]]/t^{n+1})$. Then $L^+ X = \varprojlim_n L_n^+ X$. Let $\pi_n \colon L^+ X \to L_n^+ X$ be the projection.

We call $Z \subset L^+ X(\mathbb{C})$ "fp-constructible" (respectively, "fp-open", "fp-irreducible", etc) if there exist $n$ and $Z_n \subset L_n^+ X(\mathbb{C})$ such that $Z = \pi_n^{-1}(Z_n)$, and $Z_n$ is constructible (respectively, "open", "irreducible", etc).

If $Z$ is fp-constructible, define $\overline{Z} = \pi_n^{-1}(\overline{Z_n})$. We also define $\mathrm{codim}_{L^+ X}(Z) = \mathrm{codim}_{L_n^+ X}(Z_n)$.

These definitions are independent of $n$ and $Z_n$.

2.4. **The Kazhdan-Lusztig map.** Let $O \in \underline{\mathcal{N}}$ be a nilpotent orbit. In [13, Section 9.1], Kazhdan and Lusztig show that the subset $O + t\mathfrak{g}(\mathcal{O})$ contains an open dense subset of regular semisimple elements with the same type $[w]$. Since $O + t\mathfrak{g}(\mathcal{O})$ is fp-irreducible, any two open subsets intersect. Therefore $[w]$ depends only on $O$.

Kazhdan and Lusztig also extend [13, Section 9.13] the above map to arbitrary parahoric flag varieties: given a parahoric subgroup $P$ and $O_P \in \underline{\mathcal{N}_{G_P}}$, the subset $O_P + P^+$ contains an open dense subset of regular semisimple elements with the same type $[w]$ (see also [8, Section 3.2.2]). We denote this conjugacy class by $\mathrm{KL}_G^P(O_P)$. When $P = G(\mathcal{O})$, we write $\mathrm{KL}_G(O)$.

**Remark 12.** *In [13, Section 9.13], Kazhdan and Lusztig conjecture that every conjugacy class in $W$ is in the image of some parahoric Kazhdan-Lusztig map. This only holds in Type $A_n$. According Theorem 7, the number of conjugacy classes in $W$ in the image of the Kazhdan-Lusztig map is equal to the number of 2-special* *representations of $W$. Outside of Type $A_n$, this is a strict subset of irreducible representations of $W$.*

**Example 13.** *Let $G = Sp(6)$. Conjugacy classes in $W(G)$ are in bijection with ordered pairs of partitions $(\alpha, \beta)$ with $|\alpha| + |\beta| = 3$. By considering root valuations, one can verify that the pair $(\alpha = (2), \beta = (1))$ is not in the image of a Kazhdan-Lusztig map. In this case, there are 10 irreducible representations of $W$ and only 9 2-special representations.*

2.5. **Lusztig's strata.** In [19], Lusztig defines strata of a reductive group. We recall the definition. Let $g \in G$, and let $g = su$ be its Jordan decomposition. Let $W_g$ be the Weyl group of $Z_G(s)^\circ$, and let $E_u$ be the Springer representation associated to $u$ for the reductive group $Z_G(s)^\circ$. Let $\sigma$ be the map

$$\sigma \colon G \ni g = su \to j_{W_g}^W E_u \in \{\text{2-special representations}\}.$$

**Definition 14** ([19])**.** *The fibers of $\sigma$ are strata of $G$.*

Each stratum is a union of conjugacy classes. Moreover, each stratum contains at most one unipotent conjugacy class. It is known that every stratum is locally closed ([6, Theorem 2.2]).

In [19, Section 6.4], Lusztig gives a conjectural definition of strata as follows. Let $g = su$ be as above, and let $P$ be a parahoric subgroup of $G(F)$ such that $G_P$ is of the same type as $Z_G(s)^\circ$. Lusztig conjectures that

$$\sigma' \colon G \ni g = su \to \mathrm{KL}_G^P(u)$$

is a well-defined map independent of the choice of $P$; and moreover that the fibers of $\sigma'$ are exactly the strata of $G$. These assertions follow from Theorem 7.

2.6. **Subsets of $L^\heartsuit \mathfrak{c}$.** We recall definitions and results from [27, Section 5].

Recall the discriminant function $\Delta \colon L\mathfrak{c} \to L\mathbb{A}^1$ (see Equation (2)). By definition, $L^\heartsuit \mathfrak{c}$ is the complement of the vanishing set of the discriminant function. Thus, it is an open subscheme of the closed subscheme $L^{++}\mathfrak{c} \subset L^+\mathfrak{c}$ of topologically nilpotent elements. This is the image of topologically nilpotent elements in $L^+\mathfrak{g}$ under $\chi$.

For $n > 0$, let $(L^+\mathfrak{c})_{\leq n} \subset L^+\mathfrak{c}$ be the locus such that $\mathrm{val}\,\Delta(\gamma) \leq n$ for all $\gamma \in (L^+\mathfrak{c})_{\leq n}$. Let $(L^\heartsuit \mathfrak{c})_{\leq n} = L^\heartsuit \mathfrak{c} \cap (L^+\mathfrak{c})_{\leq n}$. Clearly,

$$L^\heartsuit \mathfrak{c} = \bigcup_{n \geq 0} (L^\heartsuit \mathfrak{c})_{\leq n}.$$

and each $(L^\heartsuit \mathfrak{c})_{\leq n}$ is fp-open in $L^{++}\mathfrak{c}$.

We say that $Z \subset L^\heartsuit \mathfrak{c}$ is fp-closed (resp. fp-open, etc) if $Z \cap (L^\heartsuit \mathfrak{c})_{\leq n}$ is fp-closed (resp. fp-open, etc) as a subset of $L^+\mathfrak{c}$ for any $n \geq 0$. If $Z \subset L^\heartsuit \mathfrak{c}$ is fp-constructible, we define

$$\mathrm{codim}_{L^\heartsuit \mathfrak{c}}(Z) = \min_{n \geq 1} \mathrm{codim}_{(L^\heartsuit \mathfrak{c})_{\leq n}}(Z \cap (L^\heartsuit \mathfrak{c})_{\leq n}).$$

In this section, we drop the prefix (fp-) when the context is clear.

For any parahoric subgroup $P \subset G(F)$, $O_P \in \mathcal{N}_{G_P}$, let $(L^\heartsuit \mathfrak{c})_{O_P}^P = \chi(O_P + \mathrm{Lie}\,P^+)$ and $(L^\heartsuit \mathfrak{c})_{\overline{O_P}}^P = \chi(\overline{O_P} + \mathrm{Lie}\,P^+)$. When $P = G(\mathcal{O})$, we instead write $(L^\heartsuit \mathfrak{c})_O$ and $(L^\heartsuit \mathfrak{c})_{\overline{O}}$.

We now study these subsets of $L^\heartsuit \mathfrak{c}$. Yun proves that they are constructible and that $(L^\heartsuit \mathfrak{c})_{\overline{O_P}}^P$ is closed (Lemma 15, Lemma 16). In addition, their topological properties are related to the Kazhdan-Lusztig map ([27, Theorem 6.1], Proposition 18, Proposition 17).

**Lemma 15** ( [27, Lemma 5.2]). *Let $n \geq 1$. Suppose that $Z \subset \chi^{-1}(L^+\mathfrak{c})_{\leq n} \subset L^+\mathfrak{g}$ is fp-constructible. Then so is $\chi(Z) \subset (L^+\mathfrak{c})_{\leq n}$.*

**Lemma 16** ( [27, Lemma 5.3] or [8, Lemma 18]). *Let $P \subset LG(\mathbb{C})$ be a parahoric subgroup, and $O_P \subset \mathfrak{g}_P$ be a nilpotent orbit of $G_P$. Then $(L^\heartsuit \mathfrak{c})^P_{O_P} := \chi(\overline{O_P} + \operatorname{Lie} P^+)$ is fp-closed in $L^\heartsuit \mathfrak{c}$.*

The following propositions are a straightforward generalization of [27, Theorem 6.1] to parahoric subgroups. We reproduce the proofs with relevant changes for the convenience of the reader.

**Proposition 17.** *Let $P \subset LG(\mathbb{C})$ be a parahoric subgroup. Let $O_P \subset \mathfrak{g}_P$ be a nilpotent orbit and $[w] = \operatorname{KL}^P_G(O_P)$. Then $\delta_{[w]} = d_{O_P}$.*

*Proof.* Let $\gamma \in O_P + P^+$. By definition of the Kazhdan-Lusztig map, a generic choice of $\gamma$ satisfies $\gamma \in (L^\heartsuit \mathfrak{c})_{[w]}$. By [12, Proposition 8.2], we may further assume $\dim \operatorname{Fl}_\gamma = d_{O_P}$. Our goal is to show that a dense open subset of $O_P + P^+$ contains shallow elements of type $[w]$, for then we would have $\dim \operatorname{Fl}_\gamma = \delta_{[w]} = d_{O_P}$.

To this end, let
$$Z = \{(\gamma, gI) \mid \operatorname{Ad}(g^{-1})\gamma \in \mathfrak{I}^\heartsuit\} \subset \mathfrak{I}^\heartsuit \times \operatorname{Fl}.$$
For any $v \in \widetilde{W}$, Let $Z_v = Z \cap (\mathfrak{I}^\heartsuit \times IvI/I)$. Since Fl is stratified by the Schubert cells $IvI/I$, $Z_v$ forms a stratification of $Z$.

Note that
$$L^\heartsuit \mathfrak{I} = (\mathfrak{n} + tL^+\mathfrak{g}) \cap L^\heartsuit \mathfrak{g} = (\mathfrak{n}_P + (\mathfrak{n} \setminus \mathfrak{n}_P) + tL^+\mathfrak{g}) \cap L^\heartsuit \mathfrak{g} = (\mathfrak{n}_P + P^+) \cap L^\heartsuit \mathfrak{g}.$$
Then reduction modulo $P^+$ gives a cartesian diagram

(4)
$$\begin{array}{ccc} \bigcup_{v \in W_P} Z_v & \xrightarrow{\widetilde{\pi_P}} & \bigcup_{v \in W_P} S_v = \widetilde{\mathfrak{n}_P} \\ \downarrow{\widetilde{\operatorname{Spr}_P}} & & \downarrow{\operatorname{Spr}_P} \\ L^\heartsuit \mathfrak{I} & \xrightarrow{\pi_P} & \mathfrak{n}_P. \end{array}$$

Here, $S_v$ are conormal bundles to the Schubert cells of $G_P/B_P$; $\operatorname{Spr}_P$ is the Springer resolution for $G_P$.

The dimension for finite Springer fibers states that
$$\dim(G/B)_e = \frac{1}{2}(\dim \mathcal{N} - \dim \mathcal{O}) = \operatorname{codim}_\mathfrak{n}(\mathfrak{n} \cap \mathcal{O}).$$
Therefore, the preimage $\operatorname{Spr}_P^{-1}(\mathfrak{n}_P \cap O_P)$ contains a dense open subset $S'_v \subset S_v$ for some $v \in W_P$. Let $Z'_v = \widetilde{\pi_P}^{-1}(S'_v) \subset Z_v$. It is a dense open subset of $Z_v$.

According to [10, Section 0.8(a), 0.9], there is an open dense subset of $Z_v$ which consists of shallow elements of type $[w_v]$ for some $[w_v]$ depending on $v$. We claim that $\widetilde{\operatorname{Spr}_P}(Z'_v \cap Z''_v) \subset O_P + P^+$ has top dimension. Indeed, the map $\operatorname{Spr}_P$ is flat; hence $\operatorname{Spr}_P(S'_v)$ has top dimension in $O_P \cap \mathfrak{n}_P$. It follows that $\pi_P^{-1} \operatorname{Spr}_P(S'_v)$ has top dimension in $O_P + P^+$. Since $Z''_v$ is an open dense subset of $Z'_v$, the claim follows from the commutativity of the diagram (4).

By definition of the Kazhdan-Lusztig map, $O_P + P^+$ has a dense open subset containing elements of type $[w]$. Therefore $[w_v] = [w]$, and $O_P + P^+$ contains an open subset of shallow elements of type $[w]$. This is exactly what we needed to show. $\square$

**Proposition 18.** *Let $P \subset LG(\mathbb{C})$ be a parahoric subgroup. Let $O_P \subset \mathfrak{g}_P$ be a nilpotent orbit and $[w] = \mathrm{KL}_G^P(O)$. Then*
$$(L^\heartsuit \mathfrak{c})_{O_P}^P = \overline{(L^\heartsuit \mathfrak{c})_{[w]}^{sh}} \cap L^\heartsuit \mathfrak{c}$$

*and*
$$\mathrm{codim}_{L^\heartsuit \mathfrak{c}}(L^\heartsuit \mathfrak{c})_{O_P}^P = b_{O_P}.$$

*Finally, $(L^\heartsuit \mathfrak{c})_{O_P}^P$ is irreducible.*

*Proof.* First we claim that
$$\mathrm{codim}_{L^\heartsuit \mathfrak{c}}(L^\heartsuit \mathfrak{c})_{O_P}^P \leq d_{O_P}. \tag{5}$$

Let $\mathfrak{I}$ be the lie algebra of $I$, and $(L^\heartsuit \mathfrak{I})_{O_P}^P = \chi^{-1}(L^\heartsuit \mathfrak{c})_{O_P}^P$. According to [5, Corollary 7.4.4], the map
$$\chi \colon L^+ \mathfrak{I} \to L^+ \mathfrak{c}$$
is flat. Therefore, so is its restriction
$$\chi \colon L^\heartsuit \mathfrak{I} \to L^\heartsuit \mathfrak{c}$$
It follows that
$$\mathrm{codim}_{L^\heartsuit \mathfrak{c}}(L^\heartsuit \mathfrak{c})_{O_P}^P = \mathrm{codim}_{L^\heartsuit \mathfrak{I}}(L^\heartsuit \mathfrak{I})_{O_P}^P.$$
Let $\mathfrak{n}_P$ be the nilpotent subalgebra of $\mathfrak{g}_P$ spanned by the positive roots. To bound the right hand side, note that
$$L^\heartsuit \mathfrak{I} = (\mathfrak{n} + tL^+ \mathfrak{g}) \cap L^\heartsuit \mathfrak{g} = (\mathfrak{n}_P + (\mathfrak{n} \setminus \mathfrak{n}_P) + tL^+ \mathfrak{g}) \cap L^\heartsuit \mathfrak{g} = (\mathfrak{n}_P + P^+) \cap L^\heartsuit \mathfrak{g}.$$
Since $(L^\heartsuit \mathfrak{I})_{O_P}^P$ contains $(O_P \cap \mathfrak{n}_P) + P^+$, and $O_P \cap \mathfrak{n}_P$ has codimension $d_{O_P}$ in $\mathfrak{n}_P$, the inequality (5) follows.

Next, notice that $(L^\heartsuit \mathfrak{c})_{O_P,[w]}^P = (L^\heartsuit \mathfrak{c})_{[w]} \cap (L^\heartsuit \mathfrak{c})_{O_P}^P$ is dense in $(L^\heartsuit \mathfrak{c})_{O_P}^P$. This follows from [13, Section 9.2]: there exists a dense open subset $U \subset \overline{O_P} + P^+ \cap (L^\heartsuit \mathfrak{g})_{[w]}$. Then $\chi(U)$ is dense in $\chi(\overline{O_P} + P^+) = (L^\heartsuit \mathfrak{c})_{O_P}^P$ and by definition is contained in $(L^\heartsuit \mathfrak{c})_{[w]}$.

Therefore,
$$\delta_{[w]} = \mathrm{codim}_{L^\heartsuit \mathfrak{c}}(L^\heartsuit \mathfrak{c})_{[w]} \leq \mathrm{codim}_{L^\heartsuit \mathfrak{c}}(L^\heartsuit \mathfrak{c})_{O_P,[w]}^P = \mathrm{codim}_{L^\heartsuit \mathfrak{c}}(L^\heartsuit \mathfrak{c})_{O_P}^P \leq d_{O_P}.$$
Here, the first equality is [27, Lemma 4.2], and the final inequality is (5). By Proposition 17, equality must hold everywhere.

By [27, Lemma 4.2(3)], $(L^\heartsuit \mathfrak{c})_{[w]}^{sh}$ is dense in $(L^\heartsuit \mathfrak{c})_{[w]}$. Since $\mathrm{codim}_{L^\heartsuit \mathfrak{c}}(L^\heartsuit \mathfrak{c})_{[w]} = \mathrm{codim}_{L^\heartsuit \mathfrak{c}}(L^\heartsuit \mathfrak{c})_{O_P,[w]}^P$, the sets $(L^\heartsuit \mathfrak{c})_{[w]}$ and $(L^\heartsuit \mathfrak{c})_{O_P,[w]}^P$ must share a constructible dense set. It follows that
$$\overline{(L^\heartsuit \mathfrak{c})_{O_P,[w]}^P} \cap L^\heartsuit \mathfrak{c} = (L^\heartsuit \mathfrak{c})_{O_P}^P = \overline{(L^\heartsuit \mathfrak{c})_{[w]}^{sh}} \cap L^\heartsuit \mathfrak{c} = (L^\heartsuit \mathfrak{c}).$$

Finally, since $(L^\heartsuit \mathfrak{c})_{[w]}^{sh}$ is irreducible ([27, Lemma 4.2]), so is $(L^\heartsuit \mathfrak{c})_{O_P}^P$. $\square$

3. Affine Springer fibers and 2-special representations

In this section, we prove the equivalence of Theorem 4, Theorem 6 and Theorem 7.

3.1. **Affine Springer fibers.** It is known that $\mathrm{Fl}_\gamma^P$ is finite dimensional if and only if $\gamma$ is regular semisimple. In this case, we have

$$\dim \mathrm{Fl}_\gamma = \delta(\gamma)$$

(see (3)). This is the dimension formula of Bezrukavnikov [4].

Let $\gamma \in L^\heartsuit \mathfrak{c}$. Since the global dimension of $\mathbb{C}((t))$ is at most 1, the fiber $\chi^{-1}(\gamma) \subset L\mathfrak{g}$ is a single $LG$-conjugacy class. For $\gamma_1, \gamma_2 \in \chi^{-1}(\gamma)$ with $\gamma_1 = \mathrm{Ad}(g)\gamma_2$, left multiplication by $g$ gives an isomorphism $\mathrm{Fl}_{\gamma_1}^P \to \mathrm{Fl}_{\gamma_2}^P$.

Let $\kappa\colon \mathfrak{c} \to \mathfrak{g}$ be the Kostant section. This similarly extends to $\kappa\colon L^\heartsuit \mathfrak{c} \to L^\heartsuit \mathfrak{g}$. Let $\gamma' \in L\mathfrak{g}$, and set $\gamma = \chi(\gamma')$. By the preceding paragraph, the isomorphism class of $\mathrm{Fl}_{\gamma'}^P$ depends only on $\gamma$. In particular, it is isomorphic to $\mathrm{Fl}_{\kappa(\gamma)}^P$. We sometimes abuse notation and write $\mathrm{Fl}_\gamma^P$ when no confusion will occur.

Let $\gamma \in L^\heartsuit \mathfrak{c}$. The centralizer $Z(\gamma)$ is by definition a maximal torus in $G(F)$. Let $\mathbb{X}_*(Z(\gamma)) := \mathrm{Hom}(\mathbb{G}_m, Z(\gamma))$ be its $F$-rational cocharacter lattice. There is a map $\mathbb{X}_*(G_\gamma) \to Z(\gamma)$ given by $\lambda \mapsto \lambda(t)$. We denote its image $\Lambda = \Lambda_\gamma$.

The lattice $\Lambda$ acts on $\mathrm{Fl}_\gamma$. It is known that $Z(\gamma)$ acts freely on $\mathrm{Fl}_\gamma$, and that $\Lambda \setminus \mathrm{Fl}_\gamma$ is a projective $\mathbb{C}$-scheme.

3.2. **The $\widetilde{W}$ action on cohomology.** In [17, Section 5.4], Lusztig constructs an action of $W_{\mathrm{aff}}$ on $H_c^*(\mathrm{Fl}_\gamma)$. Yun extends this to an action of $\widetilde{W}$ in [25, Theorem 2.5]. We briefly recall the construction, following [25].

We have a Cartesian diagram

$$\begin{array}{ccc} \mathrm{Fl}_\gamma & \xrightarrow{\mathrm{ev}} & [\widetilde{\mathfrak{g}_P}/G_P] \\ \downarrow{\pi_P^I} & & \downarrow{\pi_{\mathfrak{g}_P}} \\ \mathrm{Fl}_\gamma^P & \xrightarrow{\mathrm{ev}_P} & [\mathfrak{g}_P/G_P] \end{array}$$

where $\widetilde{\mathfrak{g}_P} \to \mathfrak{g}_P$ is the Grothendieck simultaneous resolution for $G_P$. Applying proper base change,

$$\mathrm{ev}_P^* \pi_{\mathfrak{g}_P,!} \mathbb{C} = \pi_{P,!}^I \mathrm{ev}^* \mathbb{C} = \pi_{P,!}^I \mathbb{C}.$$

By classical Springer theory, $W_P$ acts on $\pi_{\mathfrak{g}_P,!}\mathbb{C}$, hence it acts on the right hand side too. This gives a $W_P$ action on $H_c^*(\mathrm{Fl}_\gamma) = \mathbf{R}\Gamma_c(\pi_{P,!}^I \mathbb{C})$. Lusztig shows that these actions glue to give a $W_{\mathrm{aff}}$ action on $H_c^*(\mathrm{Fl}_\gamma)$.

Yun (see [25, Section 2.4]) upgrades this to a $\widetilde{W}$ action as follows. Let $\Omega = X_*(T)/\mathbb{Z}\Phi^\vee$. The choice of standard Iwahori provides a splitting $\Omega \hookrightarrow \widetilde{W}$ and an identification $\Omega = N_{LG}(I)/I$. Then $\omega \in \Omega$ acts on $\mathrm{Fl}_\gamma$ by right multiplication, and hence on the left on $H_c^*(\mathrm{Fl}_\gamma)$. This $\Omega$ action and the $W_{\mathrm{aff}}$ action are compatible and generate the $\widetilde{W}$ action on $H_c^*(\mathrm{Fl}_\gamma)$.

Let $\Lambda$ be as in Section 3.1. The evaluation morphism $\mathrm{Fl}_\gamma^P \to [\mathfrak{g}_P/G_P]$ factors through the quotient $\Lambda \setminus \mathrm{Fl}_\gamma^P$ because $\Lambda$ inherits its action from $Z_{LG}(\gamma)$. To spell it out: $\mathrm{Ad}((hgI)^{-1})\gamma = \mathrm{Ad}((gI)^{-1})\gamma$ so the action of $Z_{LG}(\gamma)$ does not change the reduction type. In addition, since $\Omega$ acts on $\mathrm{Fl}_\gamma$ by right multiplication, it commutes with the $Z_{LG}(\gamma)$ action. Therefore the constructions of Lusztig and Yun extend to give a $\widetilde{W}$ action on $H_c^*(\Lambda \setminus \mathrm{Fl}_\gamma) = H^*(\Lambda \setminus \mathrm{Fl}_\gamma)$.

3.3. **2-special representations.** In this subsection, we prove the equivalence of Theorem 4, Theorem 6 and Theorem 7.

Let us give an overview of the proof. The main difficulty is in relating nilpotent orbits in the different reductive groups $G_P$, for parahorics $P \subset G(F)$. The cohomology of affine Springer fibers turns out to be a natural setting to study the Springer representations of the various $G_P$, because $W_{\text{aff}}$ contains every parahoric Weyl group $W_P$ generated by (affine) reflections.

However, for the purposes of studying Springer representation and $j$-induction, we need to study the copies of $W_P$ which are contained in $W$; the definition of affine Springer representations a priori does not work well with these subgruops. To bypass this issue, we restrict to the *stable* locus, which is the subspace of $H_c^*(\text{Fl}_\gamma)$ fixed by $Z(\gamma)$. According to a theorem of Yun [25, Theorem 1.2], the lattice $\mathbb{X}_*(T) \subset \widetilde{W}$ acts unipotently on this subspace; this allows us to identify the representations restricted to the two different copies of $W_P \subset W_{\text{aff}}$.

We start by proving a basic lemma about classical Springer theory.

**Lemma 19.** *Let $O \subset \mathcal{N}$ be a nilpotent orbit, and let $E_O$ be the corresponding Springer represntation. Let $d = \dim B_e$ be the dimension of the Springer fiber at any element $e \in O$. Then*

(1) *If $E_O$ appears in $H^*(B_{e'})$, then $e' \in \overline{O}$. In particular $\dim B_{e'} \geq \dim B_e$.*
(2) *If $E_O$ appears in $H^{d'}(B_{e'})$, then $d' \geq d$.*
(3) *The pullback map $H^d(B) \to H^d(B_e)$ is nonzero.*

*Proof.* Recall that the Springer sheaf is $\mu_* \mathbb{C}[\dim \mathcal{N}]$ where $\mu \colon \widetilde{\mathcal{N}} \to \mathcal{N}$ is the Springer resolution. We have a decomposition

$$\mu_* \mathbb{C}[\dim \mathcal{N}] = \bigoplus_{O \in \overline{\mathcal{N}}, \mathcal{L}} IC(O, \mathcal{L}) \otimes E_{O,\mathcal{L}}.$$

where $\mathcal{L}$ are irreducible local systems on $O$. By base change, we have

$$H^i(B_{e'}) = \bigoplus_{O \in \overline{\mathcal{N}}, \mathcal{L}} H^{i - \dim \mathcal{N}}(IC(O, \mathcal{L}))_{e'} \otimes E_{O,\mathcal{L}}.$$

(1) Note that $IC(O, \mathcal{L})$ is supported on $\overline{O}$. So if the stalk is nonzero, $e' \in O$.
(2) Note that the cohomology of $IC(O, \mathcal{L})$ is supported in cohomological range $[-\dim O, 0]$. So if $H^{d'}(B_{e'})$ contains $E_{O,\mathcal{L}}$, then
$$d' - \dim \mathcal{N} \geq -\dim O$$
so
$$d' \geq \dim \mathcal{N} - \dim O = d$$
by the dimension formula for Springer fibers.
(3) Since $H^*(\mathcal{N}, \mu_* \mathbb{C}[\dim \mathcal{N}]) = H^{*+\dim \mathcal{N}}(B)$, the pullback map $H^d(B) \to H^d(B_e)$ can be identified with the restriction map $H^{d - \dim \mathcal{N}}(\mathcal{N}, \mu_* \mathbb{C}[\dim \mathcal{N}]) \to H^{d - \dim \mathcal{N}}(\mu_* \mathbb{C}[\dim \mathcal{N}]_e)$. We will show that in fact
$$H^{d - \dim \mathcal{N}}(\mathcal{N}, IC(O, \mathbb{C})) \to H^{d - \dim \mathcal{N}}(IC(O, \mathbb{C})_e)$$
is nonzero. Recall that $d = \dim \mathcal{N} - \dim \mathcal{O}$, so there is a map $\mathbb{C}_{\overline{O}}[\dim O] \to IC(O, \mathbb{C})$ which induces an isomorphism on $H^{d - \dim N}$. Therefore, the above map is the composition
$$H^{d - \dim N}(\mathcal{N}, IC(O, \mathbb{C})) \cong H^{-\dim O}(\mathcal{N}, \mathbb{C}_{\overline{O}}[\dim O])$$
$$\to H^{-\dim O}(\mathbb{C}_{\overline{O}}[\dim O]_e) = \mathbb{C} \to H^{-\dim O}(IC(O, \mathbb{C})_e)$$

which is nonzero since $e \in O$.  □

**Lemma 20.** *Let $E$ be a 2-special representation of $W$; let $E = j^W_{W_P} E_P$ for some parahoric subgroup $P$ and some Springer representation $E_P$. Suppose that*

$$(E \colon \operatorname{Res}^{\widetilde{W}}_W H^i_c(\operatorname{Fl}_\gamma)^{st}) \neq 0.$$

*Then*

$$(E_P : \operatorname{Res}^{\widetilde{W}}_{W_P} H^i_c(\operatorname{Fl}_\gamma)^{st}) \neq 0.$$

*Proof.* Let $W'_P \subset \widetilde{W}$ be the parahoric subgroup of $\widetilde{W}$ generated by the subset of affine simple reflections corresponding to $P$. It is abstractly isomorphic to $W_P$; however it does not sit inside $W$.

By [25, Theorem 1.2], the action of $\mathbb{C}[X_*(T)]^W \subset \mathbb{C}[\widetilde{W}]$ factors through the action of $Z(\gamma)$. Hence it acts trivially on $H^{\text{top}}_c(\operatorname{Fl}_\gamma)$. Since $T \to T//W$ is finite, flat, and totally ramified over the point $1 \in T//W$, it follows that the $\mathbb{C}[X_*(T)]$-module $H^{\text{top}}_c(\operatorname{Fl}_\gamma)$ is supported over a power of the maximal ideal corresponding to $1 \in T$. Therefore the action of the lattice is unipotent. Now, let $V$ be the semisimplification of $H^{\text{top}}_c(\operatorname{Fl}_\gamma)^{st}$ as a $\widetilde{W}$-module. Since the lattice acts unipotently, the $W_P$-submodule structure and the $W'_P$ submodule structure of $V$ are the same. Since representations of finite groups are semisimple,

$$(E_P : \operatorname{Res}^{\widetilde{W}}_{W'_P} H^{\text{top}}_c(\operatorname{Fl}_\gamma)^{st}) = (E_P : V) = (E_P : \operatorname{Res}^{\widetilde{W}}_{W_P} H^{\text{top}}_c(\operatorname{Fl}_\gamma)^{st}).$$

Recall that $j^W_{W_P} E_P$ is a direct summand of $\operatorname{Ind}^W_{W_P} E_P$. By Frobenius reciprocity,

$$\operatorname{Hom}_W(\operatorname{Ind}^W_{W_P} E_P, \operatorname{Res}^{\widetilde{W}}_W H^{\text{top}}_c(\operatorname{Fl}_\gamma)^{st}) = \operatorname{Hom}_{W_P}(E_P, \operatorname{Res}^{\widetilde{W}}_W \operatorname{Res}^W_{W_P} H^{\text{top}}_c(\operatorname{Fl}_\gamma)^{st}) \neq 0.$$

This establishes the lemma.  □

**Lemma 21.** *Let $E = j^W_{W_P} E_P$ be a 2-special representation, such that $b_{E_P} = \dim \operatorname{Fl}_\gamma =: d$.*

*Then the following are equivalent:*
  (1) *The representation $E$ appears in $H^{\text{top}}_c(\operatorname{Fl}_\gamma)^{st}$*
  (2) *The orbit $O_P$ appears as a minimal $P$-reduction type for $\gamma$*
  (3) *The representation $E$ appears in $H^{\text{top}}(\Lambda \setminus \operatorname{Fl}_\gamma)^{st}$.*

*Proof.* (1) $\Rightarrow$ (2): According to Lemma 20, $(E_P : \operatorname{Res}^{\widetilde{W}}_{W_P} H^i_c(\operatorname{Fl}_\gamma)^{st}) \neq 0$. Let $O_P$ be the orbit corresponding to $E_P$ under the Springer correspondence.

Consider the map $\pi_P \colon \operatorname{Fl}_\gamma \to \operatorname{Fl}^P_\gamma$. Its fibers are finite Springer fibers. For any nilpotent orbit $O \subset \mathcal{N}_P$, consider $\operatorname{ev}^{-1}_P(O) = \operatorname{Fl}^{P,O}_\gamma \subset \operatorname{Fl}^P_\gamma$, and set $\operatorname{Fl}^O_\gamma = \pi^{-1}_P \operatorname{Fl}^{P,O}_\gamma$. Then we have a fibration

$$(6) \qquad B^P_O \to \operatorname{Fl}^O_\gamma \to \operatorname{Fl}^{O,P}_\gamma$$

where the fibers are finite Springer fibers in the Levi quotient. Here we use the fact that $\Lambda \subset Z(\gamma)$ centralizes $\gamma$ and in particular commutes with the evaluation map. Then

$$(7) \qquad \operatorname{Fl}_\gamma = \bigsqcup_{O \in \mathcal{N}_P} \operatorname{Fl}^O_\gamma.$$

Since the affine Springer action is constructed using the classical Springer action, the Leray-Serre spectral sequence

$$(8) \qquad E_2^{p,q} = H_c^p\left(\mathrm{Fl}_\gamma^{P,O}, H_c^q(B_O^P)\right) \Rightarrow H_c^{p+q}(\mathrm{Fl}_\gamma^O)$$

is $W_P$-equivariant, with $W_P$ acting trivially on the cohomology of the base. Since Springer fibers are projective, $H_c^q(B_O^P) = H^q(B_O^P)$. Since $E_P$ appears in $H_c^*(\mathrm{Fl}_\gamma)$, it must appear in $H_c^*(\mathrm{Fl}_\gamma^O)$ for some $O$. By Lemma 19, $O \subset \overline{O_P}$, so $d_O \geq d_{O_P}$. Additionally, we have

$$d_O \leq \dim \mathrm{Fl}_\gamma = b_{E_P} = d_{O_P}.$$

It follows that equality must hold, and so $O_P = O$. This also shows that $O_P$ is minimal.

(2) $\Rightarrow$ (3): Let $B_{O_P}^P \subset \mathrm{Fl}_\gamma$ be a fiber over some element $e \in O_P$ under the evaluation map. It embeds as a component of top dimension. Let $B^P \subset \mathrm{Fl}$ be the full flag variety containing $B_{O_P}^P$.

The set $\Sigma := \{\ell \in \Lambda \mid B \cap \ell.B \neq \varnothing\}$ is finite since $B$ is finite type. We may then replace $\Lambda$ be a sub-lattice of finite index, so that $B_{O_P}$ embeds into $\Lambda \setminus \mathrm{Fl}_\gamma$ without affecting the cohomological statements [26, Section 3.4.11].

Consider the inclusions $B_{O_P}^P \hookrightarrow B^P \hookrightarrow \Lambda \setminus \mathrm{Fl}$, where the map $B^P \hookrightarrow \mathrm{Fl}$ is the inclusion of a fiber of the fibration $\mathrm{Fl} \to \mathrm{Fl}^P$. It induces the composition in cohomology $H_c^d(\Lambda \setminus \mathrm{Fl}) \to H^d(B^P) \to H^d(B_{O_P}^P)$.

Consider the tautological line bundles on $\mathrm{Fl}$ indexed by $\mathbb{X}_*(T)$. The pullbacks of the Chern classes of these bundles generate $H^d(B^P)$, so the first map is surjective. The second map is nonzero by Lemma 19.

Therefore, this composition is a nonzero map of $W_P$-modules. Note that the image of the map $\varphi \colon H^d(B^P) \to H^d(B_{O_P}^P)$ is exactly $E_P$. Indeed, the $b$-invariant of every other representation appearing in the target is $> d$, while by definition the $b$-invariant of every representation appearing in the source is $\leq d$. It follows that there exists $v \in H^d(\Lambda \setminus \mathrm{Fl})$ such that $W_P.\varphi(v) = E_P$.

Therefore, $W.v \subset H_c^d(\Lambda \setminus \mathrm{Fl})$ is a $W$-submodule of $\mathrm{Ind}_{W_P}^W E_P$. By [11, Theorem 5.2.6],

$$\mathrm{Ind}_{W_P}^W E_P = jE_P \oplus E'$$

where any irreducible constituent of $E'$ satisfies $b_{E'} > b_{E_P}$.

On the other hand, we have the $W$-equivariant Leray-Serre fibration associated to $G/B \to \mathrm{Fl} \to \mathrm{Gr}$,

$$E_2^{p,q} = H_c^p(\mathrm{Gr}, H^q(G/B)) \Rightarrow H_c^{p+q}(\mathrm{Fl}),$$

from which it follows that any irreducible constituent of $H^d(\mathrm{Fl})$ has $b$-invariant at most $d = b_{E_P}$. Moreover, there is a Hoschild-Serre spectral sequence

$$E_2^{pq} = H^p(\Lambda, H^q(\mathrm{Fl})) \Rightarrow H^{p+q}(\Lambda \setminus \mathrm{Fl})$$

and therefore the $b$-invariant of any irreducible constituent of $H^d(\Lambda \setminus \mathrm{Fl})$ is also at most $d$ (use this fact for $H^i(\mathrm{Fl})$, $0 \leq i \leq d$). It follows that $W.v \subset H^d(\mathrm{Fl}) = j_{W_P}^W E_P$.

Now, consider the inclusion $B_{O_P}^P \hookrightarrow \Lambda \setminus \mathrm{Fl}_\gamma^P \hookrightarrow \Lambda \setminus \mathrm{Fl}$. The composition $H_c^d(\Lambda \setminus \mathrm{Fl}) \to H^d(\Lambda \setminus \mathrm{Fl}_\gamma) \to H^d(B_{O_P}^P)$ is the map $\varphi$ considered before. It follows that the image of $v$ in $H^d(\Lambda \setminus \mathrm{Fl}_\gamma)$ generates a $W$-representation isomorphic to $j_{W_P}^W E_P$.

Note that the action of $Z(\gamma)$ on $\Lambda \setminus \mathrm{Fl}_\gamma$ factors through the action of $\pi_0(Z(\gamma)/\Lambda)$, which is an étale group scheme over $\mathbb{C}$. We may assume that $\gamma \in O_P + P^+$.

Let $S = \{s_1, \ldots, s_n\}$ be the orbit of the unit coset $P/P^+ \in \Lambda \setminus \mathrm{Fl}_\gamma^{P,O_P}$, and $Z' = Z(\gamma) \cap P \subset Z(\gamma)$ be its stabilizer. Clearly, the image of $Z'$ in $P/P^+ = G_P$ is contained in $Z_{G_P}(e_P)$. By the definition of Springer representation, the induced action of $Z_{G_P}$ on $E_{O_P}$ in $H^d((G_P/B_P)_{e_P})$ is trivial. Now, since $\gamma$ is regular semisimple, $Z(\gamma)$ is a torus, in particular abelian. Therefore, $Z'$ acts on each $s_i H^d((G_P/B_P)_{e_P})$ and acts trivially on the copy of $E_{O_P}$ contained within.

For each $i$, the previous argument gives an element $w_i \in H^d(\mathrm{Fl}_\gamma)$ which generates a copy of $j_{W_P}^W E_P$. Moreover, the map $H_c^d(\Lambda \setminus \mathrm{Fl}) \to H^d(B_{O_P}^P)$ is $Z(\gamma)$ equivariant, so $Z'$ acts trivially on $w_i$. Let $V$ be the direct sum of these representations. The action of $Z(\gamma)$ descends to an action on $V$, and we just showed that $Z'$ acts trivially. Note now that $Z(\gamma)$ acts trivially on the diagonal copy of $j_{W_P}^W E_P$ in $V$. Therefore, $j_{W_P}^W E_P$ appears in $H^{\mathrm{top}}(\Lambda \setminus \mathrm{Fl}_\gamma)^{st}$.

(3) $\Rightarrow$ (1): We have the Hoschild-Serre spectral sequence
$$E_2^{pq} = H^p(\Lambda, H_c^q(\mathrm{Fl}_\gamma)) \Rightarrow H_c^{p+q}(\Lambda \setminus \mathrm{Fl}_\gamma).$$

We claim that $H^q(\mathrm{Fl}_\gamma)$ does not contain $E$ if $q < d$. If it did, then the same proof as in the proof above of (1) $\Rightarrow$ (2), combined with Lemma 19, shows that some orbit strictly smaller than $O_P$ appears; but then such an orbit has Springer fiber dimension greater than $d = \dim \mathrm{Fl}_\gamma$, a contradiction.

Therefore, the map $H^d(\Lambda \setminus \mathrm{Fl}_\gamma) \to H^0(\Lambda, H_c^d(\mathrm{Fl}_\gamma)) = H^d(\mathrm{Fl}_\gamma)^\Lambda$ is an isomorphism on the $E$-isotypic components. Moreover, this map is $Z(\gamma)$-invariant. By hypothesis, there exists a copy of $j_{W_P}^W E_P$ fixed by $Z(\gamma)$ in $H^d(\Lambda \setminus \mathrm{Fl}_\gamma)$. Therefore the same is true for $H^d(\mathrm{Fl}_\gamma)$.

□

Recall from Definition 5 that $\mathrm{Fl}_\gamma$ has property $P_{2-\mathrm{special}}$ if $H_c^{\mathrm{top}}(\mathrm{Fl}_\gamma)$ has a unique isomorphism class of 2-special representation with maximal $b$-invariant.

**Corollary 22.** *Assume that $[w]$ is in the image of some parahoric Kazhdan-Lusztig map. Suppose that Theorem 4 holds, and let $U_{[w]} \subset (L^\heartsuit \mathfrak{c})_{[w]}$ be an open set satisfying the condition in Theorem 4. Then for any $\gamma \in U_{[w]}$, $\mathrm{Fl}_\gamma$ has property $P_{2-special}$.*

*Proof.* Suppose otherwise. If $E_1 = j_{W_1}^W E_{W_1}$ and $E_2 = j_{W_2}^W E_{W_2}$ are two non-isomorphic 2-special representations with maximal $b$-invariant, then we can find parahorics $P_1$ and $P_2$, such that their Levi quotients have Weyl group $W_1$ and $W_2$ respectively. Then by Lemma 21, $\mathrm{Fl}_\gamma$ has two minimal reduction types with different $j$-induced representations. This contradicts Theorem 4. □

By running the proof of Corollary 22 backwards, it is easy to see

**Corollary 23.** *Let $U_{[w]} \subset (L^\heartsuit \mathfrak{c})_{[w]}$ be an open subset such that for any $\gamma \in U_{[w]}$, $\mathrm{Fl}_\gamma$ has property $P_{2-special}$. Then $U_{[w]}$ satisfies the requirements of Theorem 4.*

**Lemma 24.** *Suppose that $j_{W_P}^W E_P = j_{W_Q}^W E_Q$. Then $\mathrm{KL}_G^P(O_P) = \mathrm{KL}_G^Q(O_Q)$.*

*Proof.* Let $[w] = \mathrm{KL}_G^P(O_P)$. By Proposition 18, there exists a dense open subset $U_{[w]} \subset (L^\heartsuit \mathfrak{c})_{[w]}^{sh}$ such that $O_P \in \mathrm{ev}_P(\gamma)$ for any $\gamma \in U_{[w]}$. By Lemma 21, $j_{W_P}^W E_P = j_{W_Q}^W E_Q$ appears in $H^*(\mathrm{Fl}_\gamma)^{st}$. By the same lemma, $O_Q \in \mathrm{ev}_Q(\gamma)$ is a reduction type. It is necessarily minimal because

(9)  $$\delta_{[w]} = \dim \mathrm{Fl}_\gamma = d_{O_P} = b_{E_P} = b_{E_Q} = d_{O_Q}$$

by the definition of $j$-induction. This implies that $U_{[w]} \subset (L^\heartsuit \mathfrak{c})^Q_{O_Q}$. We have

$$\overline{U_{[w]}} \cap L^\heartsuit \mathfrak{c} = \overline{(L^\heartsuit \mathfrak{c})_{[w]}} \cap L^\heartsuit \mathfrak{c} = (L^\heartsuit \mathfrak{c})^P_{O_P} \subset (L^\heartsuit \mathfrak{c})^Q_{O_Q}.$$

However, by (9) and Proposition 18, the latter two sets are irreducible and have the same codimension, hence the inclusion is an equality. We now have

$$\overline{(L^\heartsuit \mathfrak{c})_{[w]}} \cap L^\heartsuit \mathfrak{c} = (L^\heartsuit \mathfrak{c})^Q_{O_Q} = \overline{(L^\heartsuit \mathfrak{c})_{\mathrm{KL}^Q_G(O_Q)}} \cap L^\heartsuit \mathfrak{c}$$

hence $[w] = \mathrm{KL}^Q_G(O_Q)$. □

**Proposition 25.** *Theorem 4 is equivalent to the following statement:*

$$\mathrm{KL}^P_G(O_P) = \mathrm{KL}^Q_G(O_Q)$$

*if and only if*

$$j^W_{W_P} E_P = j^W_{W_Q} E_Q.$$

*Proof.* First assume that Theorem 4 holds. Suppose that $\mathrm{KL}^P_G(O_P) = \mathrm{KL}^Q_G(O_Q) = [w]$. Then according to Proposition 18, we have

$$\overline{(L^\heartsuit \mathfrak{c})^{sh}_{[w]}} \cap L^\heartsuit \mathfrak{c} = (L^\heartsuit \mathfrak{c})^Q_{O_Q} = (L^\heartsuit \mathfrak{c})^P_{O_P}.$$

Let $U_{[w]} \subset (L^\heartsuit \mathfrak{c})^{sh}_{[w]}$ be an open set with the properties in Theorem 4. Since $\delta_{[w]} = \operatorname{codim}_{L^\heartsuit \mathfrak{c}} (L^\heartsuit \mathfrak{c})^P_{O_P} = \operatorname{codim}_{L^\heartsuit \mathfrak{c}} (L^\heartsuit \mathfrak{c})^Q_{O_Q}$, we may shrink $U_{[w]}$ and take it to be contained in $(L^\heartsuit \mathfrak{c})^P_{O_P} \cap (L^\heartsuit \mathfrak{c})^Q_{O_Q}$ (see [2, Lemma 2.1]). This means that for any $\gamma \in U_{[w]}$, $O_P \in \mathrm{ev}_P(\gamma)$ and $O_Q \in \mathrm{ev}_Q(\gamma)$. Moreover, since $\delta_{[w]} = d_{O_P} = d_{O_Q}$ (Proposition 17), Theorem 4 asserts that $j^W_{W_P} E_P = j^W_{W_Q} E_Q$.

The converse direction is simply Lemma 24.

Next, assume that the condition in Proposition 25 holds. Let $[w]$ be in the image of the Kazhdan-Lusztig map, say $[w] = \mathrm{KL}^P_G(O_P)$. Let $Q \subset G(F)$ be some parahoric subgroup, and let $O_Q \subset G_Q$ be a nilpotent orbit satisfying $d_{O_Q} \geq \delta_{[w]}$. Consider the intersection

$$V_{O_Q} = (L^\heartsuit \mathfrak{c})_{[w]} \cap (L^\heartsuit \mathfrak{c})^Q_{O_Q}.$$

Then $V_{O_Q}$ is a closed subset of $(L^\heartsuit \mathfrak{c})_{[w]}$. If $d_{O_Q} > \delta_{[w]}$, then it is in fact a proper closed subset (Proposition 18). If $V_{O_Q} = (L^\heartsuit \mathfrak{c})_{[w]}$, then we must have $\mathrm{KL}^Q_G(O_Q) = [w]$.

Let $S = \{(Q, O_Q) \mid d_{O_Q} >= \delta_{[w]}, V_{O_Q} \subsetneq (L^\heartsuit \mathfrak{c})_{[w]}\}$. Let

$$U_{[w]} = (L^\heartsuit \mathfrak{c})_{[w]} \setminus \bigcup_{(Q,O_Q) \in S} V_{O_Q}.$$

Then $U_{[w]}$ is a nonempty open subset of $(L^\heartsuit \mathfrak{c})_{[w]}$. Let $\gamma \in U_{[w]}$, and suppose $O_Q \in \mathrm{ev}_Q(\gamma)$ is minimal. Then $d_{O_Q} \leq \delta_{[w]} = d_{O_P}$ by choice of $U_{[w]}$. We only need to hand the reduction types with $d_{O_Q} = d_{O_P}$. In this case, by the paragraph above, $\mathrm{KL}^Q_G(O_Q) = [w] = \mathrm{KL}^P_G(O_P)$. Therefore by assumption, $j^W_{W_P} E_P = j^W_{W_Q} E_Q$. So $U_{[w]}$ satisfies the requirements of Theorem 4 □

## 4. Reductions

In this section, we make several reductions to make casework easier. To summarize, we show that to prove Theorem 4, it suffices to consider:

(1) $G$ almost simple, and for one isomorphism class of $G$ in its isogeny class
(2) only one of $G$ or $G^\vee$
(3) $P \subset G(F)$ a maximal, non-hyperspecial parahoric subgroup
(4) $[w]$ an elliptic conjugacy class of $W$.

**Lemma 26.** (1) *Let $G$ and $G'$ be isogenous reductive groups. Then for any $\gamma \in L^\heartsuit \mathfrak{c}$, $\mathrm{ev}_P^G(\gamma) = \mathrm{ev}_{P'}^{G'}(\gamma)$.*
(2) *Theorem 4 holds for $G$ if and only if it holds for every simple factor of $G$.*

*Proof.* (1) By [21, Prop. 6.6], $(\mathrm{Fl}_{G,\gamma}^P)^\circ = (\mathrm{Fl}_{G',\gamma}^{P'})^\circ$. Moreover, $Z_G(\gamma)$ and $Z_{G'}(\gamma)$ act transitively on the set of connected components. Since the action of the centralizer leaves the reduction type unchanged, it follows that $\mathrm{ev}_P^G(\gamma) = \mathrm{ev}_{P'}^{G'}(\gamma)$.

(2) $G$ is isogenous to a product of tori and simple groups. Hence this follows from (1). □

**Lemma 27.** *Let $P \subset G(F)$ be a hyperspecial parahoric subgroup of $G$. The poset $\mathrm{ev}_P(\gamma)$ is isomorphic to the poset $\mathrm{ev}_G(\gamma)$.*

*Proof.* Since $P$ is hyperspecial, the Levi quotient is isomorphic to $G$: $G_P = P/P^+ \cong G$. The affine Weyl group acts transitively on the hyperspecial vertices in the apartment. Therefore $LG_{ad}$ acts transitvely on the set of hyperspecial parahorics. Therefore the reduction type poset for any hyperspecial parahoric is isomorphic to the reduction type poset for the parahoric $P = G(\mathcal{O})$. □

**Lemma 28.** *To prove Theorem 4, it suffices to consider maximal parahoric subgroups.*

*Proof.* Let $P \subset G(F)$ be a parahoric subgroup, and $O_P$ be a minimal reduction type. Let $Q \supset P$ be a maximal parahoric subgroup containing $P$. Then $G_P$ is a Levi subgroup of $G_Q$, and we claim that the orbit $O_Q = \mathrm{Ind}_{G_P}^{G_Q} O_P$ is a reduction type for $Q$. Here Ind is the Lusztig-Spaltenstein induction [20, Theorem 1.3]. The orbit $O_Q$ satisfies $E_Q = j_{W_P}^{W_Q} E_P$ and $b_{O_Q} = b_{O_P}$ by [20, Theorem 1.3(b)]. The claim follows from these facts, and two applications of Lemma 21. Finally, hence by the transitivity of $j$ induction, $j_{W_Q}^W E_Q = j_{W_P}^W E_P$. □

Fix a $W$-equivariant isomorphism $\mathfrak{t} \cong \mathfrak{t}^*$. This allows us to identify the Chevalley base of $G$ and $G^\vee$:

$$\mathfrak{g}//G \cong \mathfrak{t}//W \cong \mathfrak{t}^*//W \cong \mathfrak{g}^\vee//G^\vee = \mathfrak{c}.$$

Under this identification, we have $\dim \mathrm{Fl}_\gamma = \dim \mathrm{Fl}_{\gamma, G^\vee}$ ([8, Lemma 20]), and $(L^\heartsuit \mathfrak{c})_{[w]}^\vee = (L^\heartsuit \mathfrak{c})_{[w]}$ ([8, Lemma 21]).

**Lemma 29.** *If $[w]$ is in the image of some parahoric Kazhdan-Lusztig map for $G$, then it is also for $G^\vee$.*

*Proof.* By assumption, we may write $[w] = \mathrm{KL}_G^P(O_P)$. Let $W(G_P) = W_x$. Then there exists an endoscopic group of $G_P$ with Weyl group $W_{x,\zeta}$, and a special representation $E_{x,\zeta}$ of $W_{x,\zeta}$ such that $E_{O_P} = j_{W_{x,\zeta}}^{W_x} E_{x,\zeta}$. Moreover, $j_{W_{x,\zeta}}^{W_\zeta} E_{x,\zeta}$ is a Springer representation for a pseudo-Levi subgroup $G_\zeta^\vee \subset G^\vee$. Let $Q_\zeta \subset LG^\vee$ be a parahoric subgroup with Levi quotient isomorphic to $G_\zeta^\vee$. Then by [8, Theorem 35], we have

$$\mathrm{KL}_{G^\vee}^{Q_\zeta}\left(\mathrm{Spr}_{G_\zeta^\vee} j_{W_{x,\zeta}}^{W_\zeta} E_{x,\zeta}\right) = \mathrm{KL}_G^{P_x}\left(\mathrm{Spr}_{G_x} j_{W_{x,\zeta}}^{W_x} E_{x,\zeta}\right) = [w],$$

so $[w]$ is also in the image of some parahoric Kazhdan-Lusztig map for $G^\vee$. □

**Proposition 30.** *Theorem 4 holds for $G$ if and only if it holds for $G^\vee$.*

*Proof.* Assume that Theorem 4 holds for $G$. We will deduce the theorem for $G^\vee$. Let $[w]$ be in the image of some parahoric Kazhdan-Lusztig map for $G$. By Lemma 29, it is in the image for $G^\vee$ as well. Let $U_{[w]} \subset (L^\heartsuit \mathfrak{c})_{[w]}$ be the open set given by Theorem 4 for $G$. By Corollary 23, it suffices to show that for any $\gamma \in U_{[w]}$, $\mathrm{Fl}_{\gamma,G^\vee}$ has property $P_{2-\mathrm{special}}$.

By Corollary 22, for any $\gamma \in U_{[w]}$, $\mathrm{Fl}_\gamma$ has property $P_{2-\mathrm{special}}$. Let $E_1$ and $E_2$ be 2-special (with respect to $G^\vee$) representations with maximal $b$-invariant, which appear in $H^{\mathrm{top}}(\Lambda \setminus \mathrm{Fl}_{\gamma,G^\vee})^{st}$ (use Lemma 21). Note that they are 2-special for $G$ as well (see Section 2.1.3). Since $[w]$ is in the image of the Kazhdan-Lusztig map for $G^\vee$,

(10) $$b_{E_1} = b_{E_2} = \dim \mathrm{Fl}_{\gamma,G^\vee} = \delta_{[w]} = \dim \mathrm{Fl}_\gamma .$$

On the other hand, by [8, Theorem 14, Prop. 31], $E_1$ and $E_2$ also appear in $H^{\mathrm{top}}(\Lambda \setminus \mathrm{Fl}_\gamma)^{st}$. Moreover, by (10) and Lemma 21, they appear in $H^{\mathrm{top}}(\mathrm{Fl}_\gamma)^{st}$, and have maximal $b$-invariant among all 2-special representations. Since $\mathrm{Fl}_\gamma$ has property $P_{2-\mathrm{special}}$, $E_1 = E_2$. This means that $\mathrm{Fl}_{\gamma,G^\vee}$ has property $P_{2-\mathrm{special}}$. □

**Lemma 31.** *Let $L$ be a Levi subgroup of $G$. If $g, g' \in L$ are in the same stratum of $L$, then they are in the same stratum of $G$.*

*Proof.* According to [19, Section 2.2], it is enough to show that

$$\mathrm{Im}(H_*(B_g) \to H_*(G/B)) = \mathrm{Im}(H_*(B_{g'}) \to H_*(G/B)).$$

We have

$$H_*(G/B) = \bigoplus_{w \in W_L \setminus W} H_*(PwB/B)$$

Let $B_g^L$ denote the Springer fiber for $L$. Note that $PwB/B \cap B_g \to B_g^L$ is an affine space bundle. Therefore, their homology can be identified up to a shift. The image $\mathrm{Im}(H_*(B_g) \to H_*(G/B))$ can be identified with the direct sum of $\mathrm{Im}(H_*((L/L \cap {}^w B)_g) \to H_*(L/(L \cap {}^w B))) = \mathrm{Im}(H_*(B_g^L) \to H_*(B^L))$. Since $g$ and $g'$ are in the same $L$-stratum, these images are equal for $g$ and $g'$, which was what we needed. □

**Proposition 32.** *Suppose that Theorem 4 holds for proper Levi subgroups of $G$. Suppose $w \in W_M$ where $M$ is a proper Levi subgroup of $G$. Then Theorem 4 holds for $[w]$.*

*Proof.* Let $Z = Z(M)^\circ$. Then $Z$ acts on $\mathrm{Fl}^P$, preserving $\mathrm{Fl}_\gamma^P$. We first prove that the $Z$-fixed points are $LM/(P \cap LM) \hookrightarrow \mathrm{Fl}^P$.

Recall from [28, Lemma 1.3.8] that we have a decomposition $LG = LB.L^+G$. If we let $B$ denote the opposite Borel to our usual choice of Borel subgroup (so $B$ contains the negative roots here), then it follows that $LG = LB.P$ where $P \subset G(F)$ is a maximal parahoric subgroup. Let $Q \supset B$ be a parabolic of $G$ containing $M$. Write $Q = M.U_Q$, then it follows that $LQ \twoheadrightarrow LG/P = \mathrm{Fl}^P$ is a surjection.

Let $[g] \in \mathrm{Fl}^P$, and suppose that it is the image of $g \in LQ$. Then $[g]$ is a fixed point if and only if $g^{-1}tg \in P \cap LQ$ for any $z \in Z$. Write $g = g'g''$, $g' \in LM$, $g'' \in LU_Q$. Let $g'' = g_1g_2$, $g_2 \in LU_Q \cap P$, $g_1 \notin P \cap LU_Q$. Explicitly, $g_2$ is a product of root groups $U_\alpha$ with some valuation constraints.

Then $g^{-1}zg \in P \cap LQ$ if and only if $g_1^{-1}zg_1 \in P \cap LB$. Let $P = P_x$ for some $x$ in the compact spherical alcove. Then $g_1$ is a product of $U_\alpha(t^{-1}\mathbb{C}[t^{-1}])$ if $\alpha(x) = 0$, and $U_\alpha(\mathbb{C}[t^{-1}])$ if $\alpha(x) < 0$, for roots $\alpha$ that are not roots of $M$. Note that the root groups in $P \cap LB$ are $U_\alpha(\mathbb{C}[[t]])$ if $\alpha(x) = 0$, and $U_\alpha(t\mathbb{C}[t])$ if $\alpha(x) < 0$. Hence $g_1^{-1}zg_1 \in P \cap LB$ if and only if $g_1^{-1}zg_1 = z$ for all $z \in Z$. It follows that $g_1 = 1$. Therefore, $[g]$ is represented by an element $g \in LM.P$. This was what we wanted.

We may assume that $M$ is chosen such that $[w]_M$ is elliptic for $M$. In particular, $[w]_M$ is in the image of the Kazhdan-Lusztig map ([27, Theorem 1.14]). Let $\gamma \in L^\heartsuit \mathfrak{t}_w$, where $LT_w \subset LM$ is a loop torus of type $w$.

Now let $\mathbb{G}_m \to Z$ be a generic cocharacter. Every point $x$ in $\mathrm{Fl}^P_\gamma$ contracts under this $\mathbb{G}_m$ action to a point $x'$ in $\mathrm{Fl}^P_{\gamma,M}$. It follows that in $\mathcal{N}_G$, the image of the evaluation map of $x'$ lies in the closure of that on $x'$. Therefore, if $O_P$ is a minimal reduction type for $P$, we may assume that it is in the image of the inclusion map from $G_P \cap M \to G_P$.

Let $(P, O_P)$ be relevant to $R(\gamma)$, i.e. $b_{O_P} = b^*(\mathrm{Fl}_\gamma)$. Since $[w]$ is in the image of the Kazhdan-Lusztig map, and shallow elements are generic in $(L^\heartsuit \mathfrak{c})_{[w]}$ we may assume $b^*(\mathrm{Fl}_\gamma) = \delta_{[w]}$. Then $\mathrm{ev}_P^{-1}(O_P) = \mathrm{Fl}_\gamma^{P,O_P}$ has dimension 0. By the above paragraph, we may assume that $O_P$ is the inclusion of some orbit $O_P^M \subset M \cap G_P$. Then $\mathrm{Fl}_{\gamma,M}^{P \cap M, O_P^M} \subset \mathrm{Fl}_\gamma^{P,O_P}$ must also be discrete. Therefore, $b_{O_P^M,M} = \dim \mathrm{Fl}_{\gamma,M}$, so it is relevant to $R_M(\gamma)$. We conclude: $O_P$ is included from a relevant orbit of $M$. If $(Q, O_Q)$ is relevant to $R(\gamma)$ as well, we need to show that $j^W_{W_Q} E_{O_Q} = j^W_{W_P} E_{O_P}$. By this same argument, $O_Q$ is the inclusion of some relevant orbit $O_Q^M$.

Now we can pick $s, s' \in T$ such that $G_Q = Z_G(s)^\circ$, $G_P = Z_G(s')^\circ$. Since the theorem holds for $M$, $sO_Q^M$ and $s'O_P^M$ satisfy $j^{W_M}_{W_{s,M}} E_{O_Q^M} = j^{W_M}_{W_{s',M}} E_{O_P^M}$. In other words, they are in the same stratum of $M$. By Lemma 31, they are in the same stratum of $G$ as well. This means that $j^W_{W_P} E_{O_P} = j^W_{W_{s'}} E_{O_P} = j^W_{W_s} E_{O_Q} = j^W_{W_Q} E_{O_Q}$ which is exactly what we needed. □

**Example 33.** *Let $G = Sp(12)$. Let $P$ be the standard parahoric subgroup with Levi quotient of type $C_3 \times C_3$. Consider the nilpotent orbits with Jordan type $O_1 = (4,2) \times (4,2)$, and $O_2 = (3,3) \times (6)$. One can check that the j-inductions of their Springer representation from $W_P$ to $W$ are equal; therefore, $\mathrm{KL}^P_G(O_1) = \mathrm{KL}^P_G(O_2)$. In particular, $\mathrm{KL}^P_G$ may not be injective, unlike $\mathrm{KL}_G$.*

*Equivalently, let $[w] = \mathrm{KL}^P_G(O_1)$. Then there is no dense open subset $U_{[w]} \subset (L^\heartsuit \mathfrak{c})_{[w]}$ such that $\mathrm{ev}_P(\gamma)$ has a unique minimal element for all $\gamma \in U_{[w]}$ (compare with Theorem 1); for a generic $\gamma \in (L^\heartsuit \mathfrak{c})_{[w]}$, both $O_1$ and $O_2$ both appear in $\mathrm{ev}_P(\gamma)$ as minimal but incomparable orbits.*

## 5. Proofs

In this section, we prove Theorem 4 for the remaining cases (Section 5.2 onwards). We first review the *skeleta* defined by Yun in Section 5.1. This is our main tool for proving Theorem 4 in the classical types. In the exceptional types, the proof relies on Alvis' tables for calculating the $j$ induction [1], as well as Spaltenstein's tables of nilpotent orbits and orbit closure diagrams [23] (among other things, see Section 5.6). In type $E_8$, we use the topology of subsets of $L^\heartsuit \mathfrak{c}$ (Proposition 18) to handle various "exceptional" cases.

### 5.1. Skeleta.

Let $w \in W$. A torus $T_w$ of type $w$ can be constructed as follows. Let $m$ be the order of $w$, and let $\zeta$ be an $m$th root of unity. Then let

$$T_w = \left( R_{\mathbb{C}((t^{1/m}))/\mathbb{C}((t))} \left( T \times_{\mathbb{C}} \mathbb{C}((t^{1/m})) \right) \right)^{\langle w \rangle, \circ}$$

where $R$ denotes the Weil restriction, $\circ$ denotes the fiberwise neutral component, and $w \in W$ acts on $T$ and multiplies $t^{1/m}$ by $\zeta$. Its Bruhat-Tits group scheme is then

$$\mathcal{T}_w = \left( R_{\mathbb{C}[[t^{1/m}]]/\mathcal{O}} \left( T \times_{\mathbb{C}} \mathbb{C}[[t^{1/m}]] \right) \right)^{\langle w \rangle}$$

Let $L^+ T_w$ be its parahoric subgroup. The following definition is due to Yun.

**Definition 34** ([27, Definition 7.1])**.** *The skeleton of type $[w]$ is the reduced structure of the fixed point subscheme of $\mathrm{Fl}$ under left translation by $L^+ T_w$*

$$X^w := (\mathrm{Fl})^{L^+ T_w, \mathrm{red}}.$$

In [27, Section 7.1] studies the skeleton in the affine Grassmanian $(\mathrm{Gr})^{L^+ T_w, \mathrm{red}}$ to understand $G(\mathcal{O})$-minimal reduction types. In this paper we study the skeleton in various parahoric partial flag varieties ($X_P^w := (\mathrm{Fl}^P)^{L^+ T_w, \mathrm{red}}$) to understand $P$-minimal reduction types.

**Lemma 35.** *For any parahoric subgroup $P \subset G(F)$,*

$$X^w \to X_P^w$$

*is surjective.*

*Proof.* If a point in $\mathrm{Fl}$ is fixed by $L^+ T_w$, then its image in $\mathrm{Fl}^P$ is certainly fixed. Let $[x] \in X_P^w$. Then $L^+ T_w$ acts on its fiber $\pi_P^{-1}([x]) \cong G_P/B_P$. Since its fiber is a flag variety (for $G_P$), it is projective. Hence the action of the connected solvable group $L^+ T_w$ has a fixed point. $\square$

**Lemma 36.** *Let $\gamma \in L^\heartsuit \mathfrak{t}_w$.*

(1) $X_P^w \subset \mathrm{Fl}_\gamma^P$.
(2) *Let $O \in \mathrm{ev}_P(\gamma)$ be minimal. Then $\mathrm{ev}_P^{-1}(O) \subset \mathrm{Fl}_\gamma$ intersects $X_P^w$.*

*Proof.*  (1) We have

$$X_P^w(\mathbb{C}) = \{gP \mid \mathrm{Ad}(g^{-1}) L^+ T_w \subset P\}.$$

(2) If $O$ is minimal, then $\mathrm{ev}_P^{-1}(O) = \mathrm{ev}_P^{-1}(\overline{O})$ is closed. Since $\mathrm{Fl}_\gamma^P$ is ind-projective, so is $\mathrm{ev}_P^{-1}(\overline{O})$. So the action of the connected solvable group $L^+ T_w$ has a fixed point.

$\square$

5.2. **Proof in Type A.** In type $A$, every maximal parahoric is hyperspecial. Moreover, every representation of $W$ is special, in particular Springer and 2-special. By Lemma 27, the minimal reduction type is the same as the minimal reduction type for $G(O)$. In particular, their Springer representations are equal.

5.3. **Proof in Type D.** Let $G = SO(V)$, where $\dim_{\mathbb{C}} V = 2n$ and $V$ is equipped with a nondegenerate quadratic form. We extend it $F$-linearly to a bilinear form $\langle -, - \rangle \colon (F \otimes V) \times (F \otimes V) \to F$. Let $V_F = V \otimes F$.

The strategy of the proof is as follows. We describe the partial flag variety $\mathrm{Fl}^P$ for various parahoric subgroups $P$ in Section 5.3.1. Next, we compute the skeleta $X_P^w$ for various elliptic elliptic elements in $W$ in Section 5.3.2. The skeleta $X_P^w$ turns out to have many points, many of which give orbits that are not minimal under the evaluation map $\mathrm{ev}_P \colon X_P^w \to [\mathfrak{g}_P/G_p]$. We identify the minimal orbits that appear in Lemma 39. Finally we prove that the $j$-induction of the Springer representations associated to these minimal orbits argree in Section 5.3.3.

5.3.1. *Parahoric flag varieties.* Let $\Lambda \subset V_F$ be a lattice. Its dual is defined to be

$$\Lambda^\vee = \{v \in V_F \mid \langle v, \Lambda \rangle \subset \mathcal{O}\}.$$

Let $\{e_i\}_{i=1}^{2n}$ be a basis of $V$ such that $\langle e_i, e_{j+n} \rangle = \delta_i^j$, $1 \leq i \leq n$. Let $\Lambda_0$ be the $\mathcal{O}$-lattice spanned by $\{e_i\}$. Then $\Lambda_0$ is self dual.

There are $n+1$ standard parahoric subgroups of $G(F)$. Let $0 \leq k \leq n$ be an integer. We will describe the $k$th standard parahoric subgroup. Let

$$\mathcal{L} = \left\{ \cdots \subset \Lambda_a \subset \Lambda_b \subset \Lambda_c \subset \ldots \;\middle|\; \begin{array}{l} \{a \leq b \leq c \ldots\} \subset \mathbb{Z}, \\ \{a, b, c \ldots\} \subset \{-k, k\} \pmod{2n}, \\ \dim_{\mathbb{C}} \Lambda_0/(\Lambda_0 \cap \Lambda_a) - \dim_{\mathbb{C}} \Lambda_a/(\Lambda_0 \cap \Lambda_a) = a \end{array} \right\},$$

and

$$\mathrm{Fl}_k = \{\Lambda_* \in \mathcal{L} \mid \Lambda_k^\vee = \Lambda_{-k},\ t\Lambda_* = \Lambda_{*+2n}\}.$$

Then $G(F)$ acts transitively on $P_k$. Consider the lattice chain
(11)
$$(te_1, \ldots, te_k, e_{k+1}, \ldots, e_{2n}) = \Lambda_k \subset \Lambda_{-k} = (t^{-1}e_{1+n}, \ldots, t^{-1}e_{k+n}, e_1, e_2, \ldots, e_{2n}).$$

It can be viewed as an element of $\mathrm{Fl}_k$.

Consider the standard parahoric $P_k \subset G(F)$ corresponding to the affine Dynkin diagram with the $k$th simple root removed. With our choice of inner product, take $T = \mathrm{diag}(t_1, t_2, \ldots, t_n, t_1^{-1}, \ldots, t_n^{-1}) \subset G$. A choice of simple roots is $\Sigma = \{\alpha_1 = \varepsilon_1 - \varepsilon_2, \ldots, \alpha_{n-1} = \varepsilon_{n-1} - \varepsilon_n, \alpha_n = \varepsilon_{n-1} + \varepsilon_n\}$. The parahoric subgroup $P_k$ corresponds to the point $x$ in the fundamental alcove such that $\alpha_i(x) = 0$ for all $i \neq k$, and $\alpha_k(x) = 1/m_k$ where $m_k$ is the coefficient of $\alpha_k$ in the highest root. One can then verify that $P_k$ is the set of matrices in $G(F)$ which can be written in the

form

$$P_k = \left(\begin{array}{cc|cc|cc|cc} \mathcal{O} & \mathcal{O} & t^{-1}\mathcal{O} & \mathcal{O} \\ \hline t\mathcal{O} & \mathcal{O} & \mathcal{O} & \mathcal{O} \\ \hline t\mathcal{O} & t\mathcal{O} & \mathcal{O} & t\mathcal{O} \\ \hline t\mathcal{O} & \mathcal{O} & \mathcal{O} & \mathcal{O} \end{array}\right) \bigcap G(F),$$

with column widths $k, n-k, k, n-k$ and row heights $k, n-k, k, n-k$.

and its pro-unipotent radical is

$$P_k^+ = \left(\begin{array}{cc|cc|cc|cc} t\mathcal{O} & \mathcal{O} & \mathcal{O} & \mathcal{O} \\ \hline t\mathcal{O} & t\mathcal{O} & \mathcal{O} & t\mathcal{O} \\ \hline t^2\mathcal{O} & t\mathcal{O} & t\mathcal{O} & t\mathcal{O} \\ \hline t\mathcal{O} & t\mathcal{O} & \mathcal{O} & t\mathcal{O} \end{array}\right) \bigcap G(F),$$

Then it is easy to see that the stabilizer of the lattice chain (11) is $P_k$. This shows that $\mathrm{Fl}_k$ can be identified with $\mathrm{Fl}/P_k$. Moreover, we have $P_k/P_k^+ \cong SO(2k) \times SO(2n-2k)$, and the action of $\gamma \in P_k$ on the Levi quotient is simply the product of its actions on $\Lambda_{-k}/\Lambda_k$ and $\Lambda_k/\Lambda_{2n-k}$.

5.3.2. *Skeleta.* Let $[w]$ be an elliptic conjugacy class in $W(D_n)$. Concretely, $[w]$ corresponds to a partition of $n$ (it represents a totally negative cycle). Then an elliptic maximal torus can be described as follows: Let $\prod_{i \in I} F_i$ be a product of extensions of $F = \mathbb{C}((t))$ of total degree $n$. For each $F_i$, let $E_i/F_i$ be a quadratic extension. Then $T_w(F) \cong \prod_i E_i^{\mathrm{Nm}=1}$, where $\mathrm{Nm}: E_i \to F_i$ is the field norm. We equip $E_i$ with the quadratic form

$$q_i(x) = \mathrm{Tr}_{F_i/F}(d_i^{-1}\mathrm{Nm}_{E_i/F_i}(x)).$$

Here $d_i$ is a generator of the different ideal of $F_i/F$. This gives a quadratic form on $V_F = \oplus E_i$ which is compatible with our earlier choice of form on $V$.

Let $\sigma_i \in \mathrm{Gal}(E_i/F_i)$ be the nontrivial element. Let $\pi_i \in \mathcal{O}_{E_i}$ be a uniformizer with $\sigma_i(\pi_i) = -1$. Then we can identify

$$L^+T_w(F) = \prod_{i \in I}(1 + \pi_i\mathcal{O}_{E_i})^{\mathrm{Nm}=1}.$$

Note that each factor $(1 + \pi_i \mathcal{O}_{E_i})^{\text{Nm}=1}$ contains some element $1 + \pi_i'$ where $\pi_i'$ is a uniformizer of $\mathcal{O}_{E_i}$. Let
$$X_k^w := (\text{Fl}_k)^{L^+ T_w}$$
denote the skeleton of type $w$ for the parahoric subgroup $P_k$. Explicitly, $X_k^w$ are lattice chains which are stable under the action of $L^+ T_w$. We now describe this explicitly.

Let
$$(\Lambda_k \subset \Lambda_{-k}) \in \text{Fl}_k$$
be a lattice chain. Define integers $a_{i,\pm k}$, $b_{i,\pm k}$ by
$$\pi_i^{a_{i,\pm k}} \mathcal{O}_{E_i} = \Lambda_{\pm k} \cap \mathcal{O}_{E_i} \text{ and } \pi_i^{b_{i,\pm k}} \mathcal{O}_{E_i} = \text{proj}_{E_i}(\Lambda_{\pm k})$$
Note that under our chosen quadratic form, $\mathcal{O}_{E_i}$ is dual to $\pi_i^{-1} \mathcal{O}_{E_i}$. From the requirement that $\Lambda_k^\vee = \Lambda_{-k}$, we get that
$$b_{i,\pm k} = -a_{i,\mp k} - 1.$$
Moreover, by multiplying with $1 + \pi_i' \in L^+ T_w$, we see that
$$b_{i,\pm k} \leq a_{i,\pm k} \leq b_{i,\pm k} + 1.$$
Since
$$\bigoplus \pi_i^{a_{i,k}} \mathcal{O}_{E_i} \subset \Lambda_k \subset \Lambda_{-k} \subset \bigoplus \pi_i^{-a_{i,k}-1} \mathcal{O}_{E_i},$$
we get that $2a_{i,k} \geq -1$, so $a_{i,k} \geq 0$. Similarly,
$$\bigoplus \pi_i^{a_{i,k}+2n_i} \mathcal{O}_{E_i} \subset t\Lambda_{-k} \subset \Lambda_k \subset \bigoplus \pi_i^{-a_{i,k}-1} \mathcal{O}_{E_i}$$
so $2a_{i,k} \geq -1 - 2n_i$ or $a_{i,k} \geq -n_i$. To summarize, up to now, we have

(12)
$$\begin{array}{ccccc}
\bigoplus \pi_i^{a_{i,-k}} \mathcal{O}_{E_i} & \subset & \Lambda_{-k} & \subset & \bigoplus \pi_i^{-a_{i,k}-1} \mathcal{O}_{E_i} \\
& & \cup & & \\
\bigoplus \pi_i^{a_{i,k}} \mathcal{O}_{E_i} & \subset & \Lambda_k & \subset & \bigoplus \pi_i^{-a_{i,-k}-1} \mathcal{O}_{E_i}
\end{array}$$

with $0 \geq a_{i,-k} \geq -n_i$, $n_i \geq a_{i,k} \geq 0$, $0 \geq a_{i,k} + a_{i,-k} \geq -1$. Conversely, it is clear that lattice chains of this form are fixed by $L^+ T_w$.

We record Spaltenstein's calculation of $d_O$ for $O$ a nilpotent orbit in Type $D_n$: [23, Prop. 6.3(V)]
$$d_O = \frac{1}{4} \sum_{i \geq 1} \left( \lambda_{2i-1}^*(\lambda_{2i-1} - 2) + (\lambda_{2i}^*)^2 \right).$$

Here $O = (\lambda_1 \geq \lambda_2 \geq \dots)$ corresponds to a partition of $2n$, and $\lambda_i^* = \#\{j \mid \lambda_j \geq i\}$ is the conjugate partition. The formula says that to maximize $d_O$, we should maximize the number of parts, and make the parts as balanced as possible.

By Lemma 27, we may assume that $1 < k < n-1$. Let us arrange $n_1 \geq n_2 \geq n_3 \dots$ in non-decreasing order. Then

**Lemma 37.** *Without loss of generality, we may assume*
(1) $a_{i,-k} = -a_{i,k}$
(2) $-a_1 \geq -a_2 \geq -a_3 \dots$
(3) $n_1 + a_1 \geq n_2 + a_2 \geq n_3 + a_3 \dots$
(4) $-n_i < a_i < 0$ if $n_i > 1$.

The proof is an easy but tedious induction argument. We omit it here.

To summarize, we are considering lattice chains of the form

$$t\Lambda_{-k} \subset \bigoplus \pi_i^{a_i-1+2n_i}\mathcal{O}_{E_i} \subset \bigoplus \pi_i^{-a_i}\mathcal{O}_{E_i} \subset \Lambda_k \subset \bigoplus \pi_i^{-a_i-1}\mathcal{O}_{E_i} \subset \bigoplus \pi_i^{a_i}\mathcal{O}_{E_i} \subset \Lambda_{-k} \subset \bigoplus \pi_i^{a_i-1}\mathcal{O}_{E_i}$$

where $a_i = a_{i,-k}$ satisfy the conditions of Lemma 37. Let $\gamma \in (L^\heartsuit \mathfrak{t}_w)^{sh}$ be a shallow element. It can be written in the form $\gamma = (\pi_i)_{i \in I}$, where $\pi_i$ is a uniformizer of $\mathcal{O}_{E_i}$. Then the action of $\gamma$ on a lattice of the above form represents the nilpotent orbit

$$(-2a_i - 1 + \delta_i) \times (2n_i - 1 + 2a_i + \delta_i')$$

where $\delta_i, \delta_i' \in \{0, 1, 2\}$. We will now determine the restrictions on $\delta_i$ and $\delta_i'$.

Let

$$U = \bigoplus \pi^{a_i-1}\mathcal{O}_{E_i}/\pi_i^{a_i}\mathcal{O}_{E_i} \text{ and } V = \bigoplus \pi^{-a_i-1}\mathcal{O}_{E_i}/\pi_i^{-a_i}\mathcal{O}_{E_i}.$$

Then our quadratic form on $V_F$ induces a perfect pairing

$$U \times V \to \mathbb{C}.$$

Let

$$A = \Lambda_{-k}/\bigoplus \pi_i^{a_i}\mathcal{O}_{E_i} \subset U \text{ and } B = \Lambda_k/\bigoplus \pi_i^{-a_i}\mathcal{O}_{E_i} \subset V.$$

The restriction of the perfect pairing to $A \times B \to \mathbb{C}$ is zero because $\Lambda_{-k} = \Lambda^\vee$. Let $U_{\geq j} = \bigoplus_{i \geq j} \pi_i^{a_i-1}$. This is a filtration on $U$. Similarly, we can define $U_{\leq j}$, $V_{\leq j}$, $V_{\geq j}$ and so on. Let $A_{\geq j} = A \cap U_{\geq j}$, and let $a_j = \dim A_{\geq j}/A_{>j}$. We similarly define $b_j$, etc.

For any $j$, consider the map $A_{\geq j} \to V_{\geq j}$ given by multiplication by $\gamma^{-2a_j}$. Similarly, consider the map $B_{\geq j} \to U_{\geq j}$ given by multiplication by $\gamma^{2n_j+2a_j}/t$. Let $A_{\geq j}^\gamma$ and $B_{\geq j}^\gamma$ be their respective images. Let $r_j, s_j$ be the ranks of the restricted pairings $A_{\geq j} \times A_{\geq j}^\gamma \to \mathbb{C}$, $B_{\geq j}^\gamma \times B_{\geq j} \to \mathbb{C}$.

**Lemma 38.** (1) *For any $j$,*

$$\sum_{j' \leq j} a_{j'} + b_{j'} \geq \frac{1}{2}(\dim U_{\leq j} + \dim V_{\leq j}) = \dim U_{\leq j}.$$

(2) *Suppose equality in (1) holds for some $j < |I|$. Let $j' > j$ be minimal such that $a_{j'} + b_{j'} > 0$. Then $r_{j'} + s_{j'} > 0$.*

*Proof.* (1) Consider the space $U \oplus V$. This is a $2|I|$ dimensional space with a quadratic form induced by the perfect pairing $U \times V \to \mathbb{C}$. Then $A \oplus B$ is Lagrangian in $U \oplus V$. Therefore, its image in $U_{\leq j} \oplus V_{\leq j}$ is co-isotropic.

(2) If equality holds, then we have a splitting

$$A_{\leq j} \oplus A_{>j} \oplus B_{\leq j} \oplus B_{>j} \subset U_{\leq j} \oplus U_{>j} \oplus V_{\leq j} \oplus V_{>j}.$$

Moreover, by choice of $j'$, $A_{>j} = A_{\geq j'}$, $B_{>j} = B_{\geq j'}$ and

$$\dim A_{\geq j'} + \dim B_{\geq j'} = \frac{1}{2}(\dim U_{\geq j'} + \dim V_{\geq j'})$$

is Lagrangian. Let $k$ be maximal such that $a_{j'} = a_{j'+1} = \cdots = a_k$. Let $l$ be maximal such that $n_{j'} - a_{j'} = n_{j'+1} - a_{j'+1} = \cdots = n_l - a_l$. Without loss of generality, we may assume $k \leq l$. Take the quotient by $U_{>k}$ and $V_{>k}$ respectively. Then the image is co-isotropic, so

$$\dim A_{k \geq * \geq j'} + \dim B_{k \geq * \geq j'} \geq \frac{1}{2}(\dim U_{k \geq * \geq j'} + \dim V_{k \geq * \geq j'}).$$

We still have a perfect pairing $\dim U_{k \geq * \geq j'} \times \dim V_{k \geq * \geq j'} \to \mathbb{C}$. Suppose that $\dim A_{k \geq * \geq j'} > \frac{1}{2} \dim U_{k \geq * \geq j'}$. Then since $\gamma$ is shallow, it is a uniformizer in each component, so $\dim A^\gamma_{k \geq * \geq j'} > \frac{1}{2} \dim V_{k \geq * \geq j'}$. Therefore, the pairing $\dim A_{k \geq * \geq j'} \times \dim A^\gamma_{k \geq * \geq j'} \to \mathbb{C}$ must be nonzero. Similarly, we are done if $\dim B_{k \geq * \geq j'} > \frac{1}{2} \dim V_{k \geq * \geq j'}$.

So we may assume that we have equalities of dimension in both cases. Assume that $A^\gamma_{k \geq * \geq j'} \perp A_{k \geq * \geq j'}$. But $B_{k \geq * \geq j'}$ is also perpendicular to $A_{k \geq * \geq j'}$ of the same dimension. So $B_{k \geq * \geq j'} = A^\gamma_{k \geq * \geq j}$. But now $A^\gamma_{k \geq * \geq j} = B_{k \geq * \geq j'} \perp B^\gamma_{k \geq * \geq j'}$ is impossible because $\gamma$ is shallow: writing $\pi_i^{2n_i} = \alpha_i t + \ldots$, the coefficients $\alpha_i$ must be distinct, so $A^\gamma_{k \geq * \geq j}$ and $B^\gamma_{k \geq * \geq j'}$ are in general position.

$\square$

**Lemma 39.** *Let $p_* = (-2a_i - 1 + \delta_i)$ and $q_* = (2n_i - 1 + 2a_i + \delta'_i)$. If $p_* \times q_*$ represents an orbit with maximal b-invariant, then $\lambda = p_* + q_*$ has the following form:*

*The partition $\lambda$ of $2n$ can be broken up into consecutive blocks, where each block is of one of the following forms:*

(1) $p_k \equiv q_k \equiv 1 \pmod{2}$ and $\lambda_k = 2n_k$ or
(2) *there exists $l > 0$ such that $p_k$ is odd, $p_{k+1}, \ldots, p_{k+2l-1}$ are even, $p_{k+2l}$ is even, $q_k$ is even, $q_{k+1}, \ldots, q_{k+2l-1}$ are even, $q_{k+2l}$ is odd and $\lambda_k = 2n_k + 1$, $\lambda_{k+i} = 2n_{k+i}$ ($1 \leq i \leq 2l - 1$) and $\lambda_{k+2l} = 2n_{k+2l} - 1$.*
(3) *same as (2) but with $p$ and $q$ reversed.*

The proof follows from an easy but tedious induction using Lemma 38; we omit it here.

5.3.3. *Characters and j-induction.* The characters of $W(D_n)$ are described by unordered pairs of partitions $(\zeta, \eta)$ with $|\zeta| + |\eta| = n$, except when $\zeta = \eta$, in which case there are two corresponding characters instead of one.

Here we work with strictly increasing sequences instead of partitions, for compatibility with [18]. The bijection is given by $p = (p_1 \leq p_2 \leq \ldots) \mapsto p + (0, 1, 2, \ldots)$.

Order the $n_i$ in non-increasing order: $n_1 \geq n_2 \geq \ldots$. Recall from [27, Theorem 9.2] that $\mathrm{RT}_{\min}(\gamma) = (2n_i + \varepsilon_i)$, where

- $\varepsilon_i = 1$ if $i$ is odd and $n_i < n_{i-1}$
- $\varepsilon_i = -1$ if $i$ is even and $n_i > n_{i+1}$
- $\varepsilon_i = 0$ otherwise

where we set $n_0 = \infty$ and $n_{|I|+1} = n_{|I|+2} = 0$. We compute the Springer representation of this orbit in Lemma 40; this handles the hyperspecial parahorics by Lemma 27. We will see that this is equal to the character of $W(D_n)$ computed above (for non-hyperspecial parahorics) (Lemma 41); this completes the proof of the theorem in Type $D_n$.

**Lemma 40.** *Order the $n_i$ in non-decreasing order: $n_1 \leq n_2 \leq n_3 \ldots$. The character corresponding to $\mathrm{RT}_{\min}(\gamma)$ is*

$$\zeta_* = \left(n_i - 1 + \left\lfloor \frac{i}{2} \right\rfloor \,\middle|\, \begin{matrix} i \text{ odd}, n_i > n_{i-1}, \\ i \text{ even}, n_i = n_{i+1} \end{matrix}\right) \text{ and } \eta_* = \left(n_i + \left\lfloor \frac{i}{2} \right\rfloor \,\middle|\, \begin{matrix} i \text{ odd}, n_i = n_{i-1}, \\ i \text{ even}, n_i < n_{i+1} \end{matrix}\right).$$

The map from orbits to characters can be found [7, p. 419]. There, the output is a partition; We briefly reproduce the version which outputs an increasing sequence here for the reader's convenience:

Let $\lambda = \lambda_1 \leq \lambda_2 \leq \ldots \lambda_{2k}$ be a partition representing an orbit. Let $\lambda_i^* = \lambda_i + i - 1$. Then $\lambda_i^*$ is a strictly increasing sequence with the same number of odd terms as even terms. Let the odd terms be

$$2\zeta_1 + 1 < 2\zeta_2 + 1 < \cdots < 2\zeta_k + 1$$

and the even parts be

$$2\eta_1 < 2\eta_2 < \cdots < \eta_k.$$

Then $(\zeta, \eta)$ is the desired pair of increasing sequences.

**Lemma 41.** *Let $(p,q)$ be partitions as in Lemma 39. Let $E$ be the Springer representation of $W_k$ corresponding to the orbit $(p,q)$. Then $j_{W_k}^W E$ is the representation appearing in Lemma 40.*

*Proof.* Let $p = (p_1 \leq p_2 \leq \cdots \leq p_{2r})$ be a partition. Let $\zeta_p$ and $\eta_p$ be defined as follows: $\zeta_p$ has one part $(p_i + i)/2 - 1$ if $i \equiv p_i \pmod{2}$, and $\eta_p$ has one part $(p_i + i - 1)/2$ if $i \not\equiv p_i \pmod{2}$. Then $(\zeta_p, \eta_p)$ is the character of $W(D_n)$ corresponding to the Springer representation of $p$.

According to [18, Section 6.3], we have $j(p,q) = (\zeta, \eta)$ where

$$\zeta_i = \zeta_{p,i} + \zeta_{q,i} - (i-1) \text{ and } \eta_i = \eta_{p,i} + \eta_{q,i} - (i-1).$$

We will process $(p,q)$ one block at a time where a block is as in Lemma 39. Suppose that our block starts with $(p_a, q_a)$. We will prove by induction that if $a$ is even, the number of elements added to $\zeta$ is $a/2$ and the number of elements added $\eta$ is $a/2 - 1$; and if $a$ is odd, the number of elements added to $\zeta = (a-1)/2$, and the number of elements added to $\eta$ is $(a-1)/2$, and that moreover, the elements agree with the first elements in $\zeta^*$ and $\eta^*$ (all elements processed in increasing order) where $(\zeta^*, \eta^*)$ is as in Lemma 40. Clearly the conditions hold when $a = 1$, in which case $\zeta$ and $\eta$ are empty.

Now let $a \geq 1$. First we deal with blocks of type (1). Recall then that $p_a$ and $q_a$ are odd, and $p_a + q_a = 2n_a$. According to the formula, if $a$ is even, then we add a new part as follows:

$$\eta_{a/2} = \frac{p_a + a - 1}{2} + \frac{q_a + a - 1}{2} - \left(\frac{a}{2} - 1\right) = \frac{p_a + q_a}{2} + \frac{a}{2} = n_a + \frac{a}{2}.$$

If $n_a < n_{a+1}$, then clearly we are done. If $n_a = n_{a+1}$, then note that $n_a + \lfloor a/2 \rfloor = n_{a+1} + \lfloor (a+1)/2 \rfloor$, and so $\eta_{a/2}$ comes from the $a+1$st part of $\mathrm{RT}_{\min}(\gamma)$. A similar calculation deals with the case where $a$ is odd.

Now let us deal with blocks of type (2). Again suppose that $k$ is even, and we have a block

$$p_k \leq p_{k+1} \leq \cdots \leq p_{k+2l}$$
$$q_k \leq \cdots \leq q_{k+2l-1} \leq q_{k+2l}$$

where $p_k, q_{k+2l}$ are odd, and the rest of the terms are even; moreover $p_k + q_k = 2n_k - 1$, $p_{k+2l} + q_{k+2l} = 2n_{k+2l} + 1$, and $p_{k+*} + q_{k+*} = 2n_{k+*}$ otherwise. Then $\eta_p$ has terms coming from $p_k, p_{k+1}, p_{k+3}, \ldots, p_{k+2l-1}$, and $\eta_q$ has terms coming from $q_{k+1}, q_{k+3}, \ldots, q_{k+2l-1}, q_{2l}$; $\zeta_p$ has terms coming from $p_{k+2}, p_{k+4}, \ldots, p_{k+2l}$, $\zeta_q$ has terms coming from $q_k, q_{k+2}, \ldots, q_{k+2l-2}$. Moreover, from the description of blocks, since $p$ and $q$ are orthogonal partitions, we must have $q_{k+2i} = q_{k+2i+1}$ and $p_{k+2i+1} = p_{k+2i+2}$ for $0 \leq i < l$.

By inductive hypothesis, $\zeta$ has $k/2$ terms and $\eta$ has $k/2 - 1$ terms. We have

$$\eta_{k/2} = \frac{p_k + k - 1}{2} + \frac{q_{k+1} + (k+1) - 1}{2} - \left(\frac{k}{2} - 1\right) = \frac{p_k + q_k + 1}{2} - \frac{k}{2} = n_i - \frac{k}{2},$$

and for $(0 \leq \alpha < l - 1)$,

$$\eta_{k/2+\alpha+1} = \frac{p_{k+1+2\alpha} + k + 1 + 2\alpha - 1}{2} + \frac{q_{k+3+2\alpha} + k + 3 + 2\alpha - 1}{2} - \left(\frac{k}{2} + \alpha\right)$$
$$= \frac{p_{k+2+2\alpha} + q_{k+2+2\alpha}}{2} + \frac{k}{2} + \alpha + 1 = n_{k+2+2\alpha} + \frac{k}{2} + \alpha + 1,$$

and

$$\eta_{k/2+l} = \frac{p_{k+2l-1} + k + 2l - 1 - 1}{2} + \frac{q_{k+2l} + k + 2l - 1}{2} - \left(\frac{k}{2} + l - 1\right)$$
$$= \frac{p_{k+2l} + q_{k+2l} - 1}{2} + \frac{k}{2} + l = n_{k+2l} + \frac{k}{2} + l,$$

and finally $(0 < \alpha \leq l)$

$$\zeta_{k/2+\alpha} = \frac{p_{k+2\alpha} + k + 2\alpha}{2} - 1 + \frac{q_{k+2(\alpha-1)} + k + 2(\alpha-1)}{2} - 1 - \left(\frac{k}{2} + \alpha - 1\right)$$
$$= \frac{p_{k+2\alpha} + q_{k+2(\alpha-1)}}{2} + \frac{k}{2} + \alpha - 2$$
$$= \frac{p_{k+2\alpha-1} + q_{k+2\alpha-1}}{2} + \frac{k}{2} + \alpha - 2$$
$$= n_{k+2\alpha-1} + \left\lfloor \frac{k + 2\alpha - 1}{2} \right\rfloor - 1.$$

This establishes the inductive step for $k$ even. The case $k$ odd, and blocks of type (3) are handled similarly. This finishes the inductive step. □

5.4. **Proof in Type B.** The elliptic conjugacy classes $W(B_n)$ have the same description as in $W(D_n)$. One can show that the parahoric versions of the skeleta have similar descriptions as well. Therefore, in principle, the same proof strategy as in Section 5.3 should work.

However, the situation is now more complicated. Namely, let $[w]$ be an elliptic conjugacy class corresponding to a partition of $n$; then we can identify

$$V_F = F \oplus \bigoplus_{i \in I} E_i$$

where $E_i$ are as in Section 5.3.2. The quadratic form on $E_i$ remains the same; however, to ensure that the discriminant of the resulting form on $V_F$ has even valuation, we are forced to equip $F$ with different quadratic forms depending on whether $|I|$ is even or odd. The calculation of the skeleton is slightly different in each case.

To bypass this casework, we choose to use the compatibility of Theorem 4 with Langlands duality (Proposition 30) and deduce it from the Type C case.

5.5. **Proof in Type C.** Let $G = Sp(V)$, where $\dim_{\mathbb{C}}(V) = 2n$ and $V$ is equipped with a nondegenerate symplectic form. Let $V_F = V \otimes F$. Let $\langle -, - \rangle$ denote the induced symplectic form. The calculation in this case are almost completely parallel to the calculations in Type D. We briefly give the argument, highlighting the parts that differ.

5.5.1. *Parahoric flag varieties.* The descriptions in Type $D_n$ carry over almost directly. We record the answer for the reader's convenience:

$$\mathrm{Fl}_k = \{\Lambda_k \subset \Lambda_{-k} \mid \Lambda_k^\vee = \Lambda_{-k}\}$$

where again $\Lambda^\vee = \{v \in V_F \mid \langle v, \Lambda \rangle \subset O\}$. The Levi quotient $P_k/P_k^+$ is isomorphic to $Sp(2k) \times Sp(2(n-k))$. For any $\gamma \in P_k$, $\gamma$ induces an action on $\Lambda_k/\Lambda_{-k} \cong \mathbb{C}^{2k}$ and $\Lambda_{-k}/t\Lambda_k \cong \mathbb{C}^{2(n-k)}$. This action represents the image of $\gamma$ in the Levi quotient.

5.5.2. *Skeleta.* Let $[w]$ be an elliptic conjugacy class in $W(C_n)$. Like in the Type $D_n$ case, it represents a negative cycle with lengths $(n_1, \ldots n_r)$, $n_1 + \cdots + n_r = n$. Let $\prod_{i \in I} F_i$ be a product of extensions of $F = \mathbb{C}((t))$ of total degree $n$. For each $F_i$, let $E_i/F_i$ be a quadratic extension. Then $T_w(F) \cong \prod_i E_i^{\mathrm{Nm}=1}$, where $\mathrm{Nm} \colon E_i \to F_i$ is the field norm. We equip $E_i$ with the symplectic form

$$\langle x, y \rangle_i = \mathrm{Tr}_{F_i/F} \left( \frac{x\sigma_i(y) - y\sigma_i(x)}{\pi_i d_i} \right).$$

Here $d_i$ is a generator of the different ideal of $F_i/F$ and $\sigma_i \in \mathrm{Gal}(E_i/F_i)$ is the nontrivial element. Under this form $\mathcal{O}_{E_i}$ is a self-dual lattice. This is the main difference from the Type $D_n$ case. The same argument as in Section 5.3.2 shows

**Lemma 42.**

$$X_k^w := (\mathrm{Fl}_k)^{L^+T_w} = \left\{ \begin{array}{ccccc} \bigoplus \pi_i^{a_{i,-k}} \mathcal{O}_{E_i} & \subset & \Lambda_{-k} & \subset & \bigoplus \pi_i^{-a_{i,k}} \mathcal{O}_{E_i} \\ & & \cup & & \\ \bigoplus \pi_i^{a_{i,k}} \mathcal{O}_{E_i} & \subset & \Lambda_k & \subset & \bigoplus \pi_i^{-a_{i,-k}} \mathcal{O}_{E_i} \end{array} \right\}$$

where $0 \leq a_{i,-k} + a_{i,k} \leq +1$. The difference from the Type $D_n$ calculation is that $\mathrm{pr}_i(\Lambda_{-k}) = \pi_i^{a_{i,k}} \mathcal{O}_{E_i}$, instead of $\mathrm{pr}_i(\Lambda_{-k}) = \pi_i^{a_{i,k}-1} \mathcal{O}_{E_i}$ (compare (12)).

The formula for $d_O$ for $O = (\lambda_1 \geq \lambda_2 \geq \ldots)$ a nilpotent orbit of $Sp(V)$ is (from [23, Prop. 6.3(I)])

$$d_O = \frac{1}{4} \sum_{i \geq 1} \left( (\lambda_{2i-1}^*)^2 + \lambda_{2i}^*(\lambda_{2i}^* - 2) \right).$$

By Lemma 27, we may assume that $0 < k < n$. Let us arrange $n_1 \geq n_2 \geq n_3 \ldots$ in non-decreasing order. We have the analogue of Lemma 37:

**Lemma 43.** *Without loss of generality, we may assume*

(1) $a_{i,-k} = 1 - a_{i,k}$
(2) $-a_1 \geq -a_2 \geq -a_3 \ldots$
(3) $n_1 + a_1 \geq n_2 + a_2 \geq n_3 + a_3 \ldots$
(4) $-n_i < a_i < 0$ *if* $n_i > 1$.

Let $\gamma \in (L^\heartsuit \mathfrak{t}_w)^{sh}$ be a shallow element. It can be written in the form $\gamma = (\pi_i)_{i \in I}$, where $\pi_i$ is a uniformizer of $\mathcal{O}_{E_i}$. The action of $\gamma$ on a lattice of the above form represents the nilpotent orbit

$$(-2a_i + \delta_i) \times (2n_i - 2 + 2a_i + \delta'_i)$$

where $\delta_i, \delta'_i \in \{0, 1, 2\}$ and $a_i = a_{i,-k}$. The following lemma is the analogue of Lemma 39. The main difference is that the parities in the conditions are swapped. This is a consequence of the fact that $\mathcal{O}_{E_i}$ is self dual.

**Lemma 44.** *Let $p_* = (-2a_i + \delta_i)$ and $q_* = (2n_i + 2a_i - 2 + \delta'_i)$. If $p_* \times q_*$ represents an orbit with maximal b-invariant, then $\lambda = p_* + q_*$ has the following form:*

*The partition $\lambda$ of $2n$ can be broken up into consecutive blocks, where each block is of one of the following forms:*

  (1) $p_k \equiv q_k \equiv 0 \pmod{2}$ *and* $\lambda_k = 2n_k$ *or*
  (2) *there exists $l > 0$ such that $p_k$ is even, $p_{k+1}, \ldots, p_{k+2l-1}$ are odd, $p_{k+2l}$ is odd, $q_k$ is odd, $q_{k+1}, \ldots, q_{k+2l-1}$ are odd, $q_{k+2l}$ is even and $\lambda_k = 2n_k + 1$, $\lambda_{k+i} = 2n_{k+i}$ ($1 \leq i \leq 2l - 1$) and $\lambda_{k+2l} = 2n_{k+2l} - 1$.*
  (3) *same as (2) but with $p$ and $q$ reversed.*

5.5.3. *Characters and j-induction.* Recall from [27, Section 8.4] that $\mathrm{RT}_{\min}(\gamma) = (2n_i)$.

The characters of $W(C_n)$ are described by ordered pairs of partitions $(\zeta, \eta)$ with $|\zeta| + |\eta| = n$.

**Lemma 45.** *Order the $n_i$ in non-decreasing order: $n_1 \leq n_2 \leq n_3 \cdots \leq n_{2r}$, adding a 0 part if necessary. The character corresponding to $\mathrm{RT}_{\min}(\gamma)$ is*

$$\zeta_* = (n_{2i} + i - 1 \mid 1 \leq i \leq r) \text{ and } \eta_* = (n_{2i-1} + i - 1 \mid 1 \leq i \leq r).$$

The Springer correspondence in Type $C_n$ is the same as in Type $D_n$, except that we may need to pad $\lambda$ with a leading 0 to ensure it has an even number of parts [7, Page 419]. Lemma 45 follows easily. Finally, we prove

**Lemma 46.** *Let $(p, q)$ be partitions as in Lemma 44. Let $E$ be the Springer representation of $W_k$ corresponding to the orbit $(p, q)$. Then $j^W_{W_k} E$ is the representation appearing in Lemma 45.*

*Proof.* Let $p = (p_1 \leq p_2 \leq \cdots \leq p_{2r})$ be a partition. Let $\zeta_p$ and $\eta_p$ be defined as follows: $\zeta_p$ has one part $(p_i + i)/2 - 1$ if $i \equiv p_i \pmod{2}$, and $\eta_p$ has one part $(p_i + i - 1)/2$ if $i \not\equiv p_i \pmod{2}$. Then $(\zeta_p, \eta_p)$ is the character of $W(C_n)$ corresponding to the Springer representation of $p$.

According to [18, Section 4.5], we have $j(p, q) = (\zeta, \eta)$ where

$$\zeta_i = \zeta_{p,i} + \zeta_{q,i} - (i - 1) \text{ and } \eta_i = \eta_{p,i} + \eta_{q,i} - (i - 1).$$

We will process $(p, q)$ one block at a time where a block is as in Lemma 44. Suppose that our block starts with $(p_a, q_a)$. We will prove by induction that if $a$ is even, the number of elements added to $\zeta$ is $a/2 - 1$ and the number of elements added $\eta$ is $a/2$; and if $a$ is odd, the number of elements added to $\zeta = (a - 1)/2$, and the number of elements added to $\eta$ is $(a - 1)/2$, and that moreover, the elements agree with the first elements in $\zeta^*$ and $\eta^*$ (all elements processed in increasing order) where $(\zeta^*, \eta^*)$ is as in Lemma 40. Clearly the conditions hold when $a = 1$, in which case $\zeta$ and $\eta$ are empty.

Now let $a \geq 1$. First we deal with blocks of type (1). Recall then that $p_a$ and $q_a$ are even, and $p_a + q_a = 2n_a$. According to the formula, if $a = 2i$ is even, then we add a new part as follows:

$$\zeta_i = \frac{p_{2i} + 2i}{2} - 1 + \frac{q_{2i} + 2i}{2} - 1 - (i-1) = \frac{p_{2i} + q_{2i}}{2} - 1 + i = n_{2i} - 1 + i$$

A similar calculation deals with the case where $a$ is odd.

Now let us deal with blocks of type (2). Again suppose that $k = 2i$ is even, and we have a block

$$p_k \leq p_{k+1} \leq \cdots \leq p_{k+2l}$$
$$q_k \leq \cdots \leq q_{k+2l-1} \leq q_{k+2l}$$

where $p_k, q_{k+2l}$ are even, and the rest of the terms are odd; moreover $p_k + q_k = 2n_k - 1$, $p_{k+2l} + q_{k+2l} = 2n_{k+2l} + 1$, and $p_{k+*} + q_{k+*} = 2n_{k+*}$ otherwise. Then $\zeta_p$ has terms coming from $p_k, p_{k+1}, p_{k+3}, \ldots, p_{k+2l-1}$, and $\zeta_q$ has terms coming from $q_{k+1}, q_{k+3}, \ldots, q_{k+2l-1}, q_{2l}$; $\eta_p$ has terms coming from $p_{k+2}, p_{k+4}, \ldots, p_{k+2l}$, $\eta_q$ has terms coming from $q_k, q_{k+2}, \ldots, q_{k+2l-2}$. Moreover, from the description of blocks, since $p$ and $q$ are symplectic partitions, we must have $q_{k+2i} = q_{k+2i+1}$ and $p_{k+2i+1} = p_{k+2i+2}$ for $0 \leq i < l$.

By inductive hypothesis, $\zeta$ has $k/2 - 1$ terms and $\eta$ has $k/2$ terms. We have

$$\zeta_i = \frac{p_{2i} + 2i}{2} - 1 + \frac{q_{2i+1} + 2i + 1}{2} - 1 - (i-1) = \frac{p_{2i} + q_{2i} + 1}{2} + i - 1 = n_{2i} + i - 1,$$

and for $(0 \leq \alpha < l - 1)$,

$$\zeta_{i+\alpha+1} = \frac{p_{2i+2\alpha+1} + 2i + 2\alpha + 1}{2} - 1 + \frac{q_{2i+2(\alpha+1)+1} + 2i + 2(\alpha+1) + 1}{2} - 1 - (i+\alpha)$$
$$= \frac{p_{2i+2\alpha+2} + q_{2i+2\alpha+2}}{2} + (i+\alpha) = n_{2i+2\alpha+2} + i + \alpha,$$

and

$$\zeta_{i+l} = \frac{q_{2i+2l} + 2i + 2l}{2} - 1 + \frac{p_{2i+2l-1} + 2i + 2l - 1}{2} - 1 - (i+l-1)$$
$$= \frac{q_{2i+2l} + p_{2i+2l} - 1}{2} + i + l - 1 = n_{2i+2l} + i + l - 1,$$

and finally $(0 < \alpha \leq l)$

$$\eta_{i+\alpha+1} = \frac{q_{2i+2\alpha} + 2i + 2\alpha - 1)}{2} + \frac{p_{2i+2(\alpha+1)} + 2i + 2\alpha + 2 - 1}{2} - (i + \alpha + 1 - 1)$$
$$= \frac{q_{2i+2\alpha+1} + p_{2i+2\alpha+1}}{2} + i + \alpha$$
$$= n_{2i+2\alpha+1} + i + \alpha$$

This establishes the inductive step for $k$ even. The case $k$ odd, and blocks of type (3) are handled similarly. This finishes the inductive step. $\square$

5.6. **Proof in Exceptional types.** In this section, we prove Theorem 4 for groups of exceptional type. We first outline the strategy.

Let $G$ be of exceptional type, and let $[w] \in W(G)$ be an elliptic conjugacy class, and suppose we want to prove Theorem 4 for the class $[w]$. Then by [27, Theorem 1.14], $[w] = \mathrm{KL}_G(O_{[w]})$ for some nilpotent orbit $O_{[w]}$. Let $P \subset G(F)$ be a parahoric subgroup. We may assume that $P$ is maximal. Let $O_P \in G_P$ be a nilpotent orbit. It suffices to prove that

(13) $\quad \mathrm{KL}_G(O_{[w]}) = [w] = \mathrm{KL}_G^P(O_P)$ if and only if $j_{W_P}^W E_{O_P} = E_{O_{[w]}}$

by Proposition 25. By Proposition 17, we only need to check (13) for orbits such that $d_{O_{[w]}} = d_{O_P}$.

**Lemma 47.** *Suppose that $d_{O_{[w]}} = d_{O_P}$, and $j_{W_P}^W E_{O_P}$ is a Springer representation. Then* (13) *holds for* $([w], O_P, P)$.

*Proof.* Let $O'$ be the orbit corresponding to the Springer representation $j_{W_P}^W E_{O_P}$.

Assume that $\mathrm{KL}_G(O_{[w]}) = [w] = \mathrm{KL}_G^P(O_P)$. Then by Proposition 18, there exists an open dense subset $U \subset (L^\heartsuit \mathfrak{c})_{[w]}^{sh}$ such that for any $\gamma \in U$, $O_P \in \mathrm{ev}_P(\gamma)$ and $O_{[w]} \in \mathrm{ev}_G(\gamma)$. Applying Lemma 21 twice, we see that $O' \in \mathrm{ev}_G(\gamma)$. Since $d_{O_{[w]}} = d_{O_P}$, $O'$ and $O_{[w]}$ are both minimal elements; but according to Yun's theorem (Theorem 1, [27, Theorem 1.11]), the minimal element is unique. Therefore $O' = O_{[w]}$, and their Springer representations have to be equal.

On the other hand, assume that $j_{W_P}^W E_{O_P} = E_{O_{[w]}}$. Then by Proposition 18, there exists an open dense subset $U \subset (L^\heartsuit \mathfrak{c})_{[w]}^{sh}$ such that for any $\gamma \in U$, $O \in \mathrm{ev}_G(\gamma)$. By Lemma 21, $O_P \in \mathrm{ev}_P(\gamma)$. Therefore,

$$(L^\heartsuit \mathfrak{c})_{O_P}^P \supset \overline{U} = (L^\heartsuit \mathfrak{c})_{\overline{O_{[w]}}}.$$

By Proposition 18, these sets are irreducible of the same dimension, hence they must be equal. By definition, $\mathrm{KL}_G(O_{[w]}) = \mathrm{KL}_G^P(O_P)$. □

**Corollary 48.** *Suppose that $d_{O_{[w]}} = d_{O_P}$, and $O_P$ is special. Then* (13) *holds for* $([w], O_P, P)$.

*Proof.* According to [16, Theorem 10.7], if $O_P$ is special then $j_{W_P}^W E_{O_P}$ is a Springer representation. □

We hence give the following procedure to verify Theorem 4 in exceptional types. Fix $G$ of exceptional type.

(1) Enumerate all elliptic conjugacy classes $[w]$ of $G$, and find $\delta_{[w]}$ for each.
(2) Let $G_P$ be the Levi quotient of a maximal parahoric subgroup $P$. For any $P$, let

$$S_P = \{O_P \in \underline{\mathcal{N}_{G_P}} \mid d_{O_P} = \delta_{[w]} \text{ for some } [w], \text{ and } O_P \text{ not special}\}.$$

(3) If, for every $O_P \in S_P$, $j_{W_P}^W O_P$ is Springer, then we are done.
(4) Let $S'_P = \{O_P \in S_P \mid j_{W_P}^W O_P \text{ is not Springer}\} \subset S_P$. Check that $\mathrm{KL}_G^P(O_P) \neq [w]$ for any elliptic $[w]$ with $d_{O_P} = \delta_{[w]}$.

In fact, we will see that condition (3) holds if $G$ is not of type $E_8$.

Now we list our sources for the data required. The elliptic conjugacy classes are enumerated in [14, Section 4.3]. There, Lusztig lists $[w], \Phi([w])$ for $[w]$ and elliptic conjugacy class and $\Phi \colon \underline{W} \to \underline{\mathcal{N}}$ a map from conjugacy classes in $W$ to nilpotent orbits. In [27, Theorem 1.14], Yun shows that $\Phi([w]) = \mathrm{RT}_{\min}([w])$. Therefore, $\delta_{[w]} = d_{\Phi([w])}$ by [27, Theorem 6.1].

For nilpotent orbits in classical groups, [9, Section 6.3] gives criteria for orbits to be special. The values of $d_O$ can be deduced from [23, Section 6.3]. The Springer map from an orbit to a character of $W$ is given in [7, Chapter 13.3].

In [22], Spaltenstein gives a table of nilpotent orbits in exceptional groups, along with $d_O$ as well as the character of $W$ associated to a nilpotent orbit. In particular we can check if a character corresponds to a Springer representation. Note that the notation for characters of $W$ is slightly different from [7, Chapter 13.3]. The $b$-invariant of a character can be deduced from [3]. Using this, one can calculate the $j$-induction of characters using [1].

5.6.1. *Type $G_2$.* A maximal parahoric of $G$ has Levi quotient of type $A_1 \times A_1$, $G_2$ or $A_2$. In the type $A_n$ cases, every orbit is special. So Corollary 48 handles these cases. The parahoric with Levi quotient $G_2$ is hyperspecial; we do not need to check hyperspecial parahorics by Lemma 27. This handles type $G_2$.

5.6.2. *Type $E_6$.* A maximal parahoric of $G$ has Levi quotient of type $A_1 \times A_5$, $E_6$, $A_2 \times A_2 \times A_2$. Therefore this case is also handled by Corollary 48 and Lemma 27.

5.6.3. *Type $F_4$.* A maximal parahoric of $G$ has Levi quotient of type $A_1 \times C_3$, $B_4$, $A_2 \times A_2$, $A_3 \times A_1$. We only need to handle the first two cases (again by Corollary 48). The basic orbits (in bijection with elliptic conjugacy classes in $W$) in type $F_4$ are listed in Table 1.

| Basic Nilpotent Orbit | $F_4$ | $F_4(a_1)$ | $F_4(a_2)$ | $B_3$ | $C_3$ | $F_4(a_3)$ | $C_3(a_1)$ | $\widetilde{A_2} + A_1$ | $A_1 + \widetilde{A_1}$ |
|---|---|---|---|---|---|---|---|---|---|
| $d_O = \delta_{[w]}$ | 0 | 1 | 2 | 3 | 3 | 4 | 5 | 6 | 10 |

TABLE 1. Basic nilpotent orbits in $F_4$

First we handle the case of $A_1 \times C_3$. Every orbit in type $A_1$ is special so we only consider non-special orbits in $C_3$. Thus we only consider orbits corresponding to $(4, 1, 1)$ or $(2, 1, 1, 1, 1)$. The orbit $O = (2, 1, 1, 1, 1) \times (1, 1)$ has $d_O = 7$; we do not need to consider this orbit in step (2). The rest of the orbits have Springer $j$-induction (see [22]). We tabulate results in Table 2.

| Orbit | $d_O$ | Character | $j\chi$ |
|---|---|---|---|
| $(4, 1, 1) \times (2)$ | 2 | $((2, 1), ()) \times (2)$ | $\chi_{9,1}$ |
| $(2, 1^4) \times (2)$ | 6 | $((1^3), ()) \times (2)$ | $\chi_{9,2}$ |
| $(4, 1, 1) \times (1, 1)$ | 3 | $((2, 1), ()) \times (1, 1)$ | $\chi_{8,1}$ |

TABLE 2. Case $C_3 \times A_1$

Next we handle the case $B_4$. The non-special orbits are $(2^2 1^5)$, $(2^4 1)$, $(4^2 1)$ with $d_O = 10, 8, 3$ respectively. According Table 1, we do not need to consider the orbit $(2^4 1)$. We tabulate results in Table 3.

| Orbit | $d_O$ | Character | $j\chi$ |
|---|---|---|---|
| $2^2 1^5 1$ | 10 | $((), (2, 1, 1))$ | $\chi_{9,4}$ |
| $4^2 1$ | 3 | $((1), (3))$ | $\chi_{8,1}$ |

TABLE 3. Case $B_4$

This finishes the proof in type $F_4$.

5.6.4. *Type $E_7$.* One checks that if $P \subset G(F)$ is a maximal parahoric, then its Levi quotient is either a product of groups of type $A_n$, or it is hyperspecial, or it is of type $D_6 \times A_1$. Note that there are two embeddings of $W(D_6 \times A_1) \subset W(E_7)$; however, they differ by a graph automorphism of the affine Dynkin diagram which induces the identity function on $\mathrm{Irr}(W)$ ([1, Page 28]). So it suffices to consider only one embedding.

The basic orbits in type $E_7$ are as in Table 4.

| Basic Orbit | $\delta_{[w]} = d_O$ | Basic Orbit | $\delta_{[w]} = d_O$ |
|---|---|---|---|
| $E_7$ | 0 | $D_6(a_2) + A_1$ | 7 |
| $E_7(a_1)$ | 1 | $D_6(a_2)$ | 8 |
| $E_7(a_2)$ | 2 | $(A_5 + A_1)'' = A_5 + A_1$ | 9 |
| $D_6 + A_1 = E_7(a_3)$ | 3 | $D_4 + A_1$ | 12 |
| $D_6$ | 4 | $A_3 + A_2 + A_1$ | 13 |
| $D_6(a_1) + A_1 = E_7(a_4)$ | 5 | $4A_1$ | 28 |

Table 4. Basic orbits in $E_7$

Note that the notation for orbit names differs sometimes between [22] and [14]. We write down both names when this happens. Again, every orbit in type $A_1$ is special so we only consider non-special orbits in $D_6$. We tablulate results in Table 5.

| Orbit | $d_O$ | Character | $j\chi$ |
|---|---|---|---|
| $(1^5 2^2 3^1) \times (1^2)$ | 15 | Not a value of $\delta_{[w]}$ | - |
| $(1^5 2^2 3^1) \times (2)$ | 14 | Not a value of $\delta_{[w]}$ | - |
| $(1^3 2^2 5^1) \times (1^2)$ | 9 | $((1^1 2^1 3^1), ()) \times (1^2)$ | $\phi_{216,9} = 216_{a'}$ |
| $(1^3 2^2 5^1) \times (2)$ | 8 | $((1^1 2^1 3^1), ()) \times (2)$ | $\phi_{280,8} = 280_b$ |
| $(1^1 2^4 3^1) \times (1^2)$ | 13 | $((2^3), ()) \times (1, 1)$ | $\phi_{210,13} = 210_b \otimes \mathrm{sign}$ |
| $(1^1 2^4 3^1) \times (2)$ | 12 | $((2^3), ()) \times (2)$ | $\phi_{84,12} = 84_a$ |
| $(1^1 2^2 7) \times (1^2)$ | 5 | $((), (4, 2)) \times (1^2)$ | $\phi_{189,5} = 189'_b$ |
| $(1^1 2^2 7) \times (2)$ | 4 | $((), (4, 2)) \times (2)$ | $\phi_{35,4} = 35_b$ |

Table 5. Case $D_6 \times A_1$

Note that sometimes the notation for characters of $W$ differs between [22] and [3]. We write down both when this occurs. This finishes the proof in type $E_7$.

5.6.5. *Type $E_8$.* The relevant parahorics have Levi quotient of type $D_8$, $A_3 \times D_5$, $A_2 \times E_6$, $A_1 \times E_7$. The basic orbits in type $E_8$ are as in Table 6.

Below, there may be orbits whose $j$-induction is not Springer (see Item 4). We indicate these with ⚠ and we handle them in Section 5.6.6.

In type $D_8$, the special orbits have $d_O$ equal to $4, 7, 8, 10, 11, 12, 13, 18, 19, 26, 14, 14, 22, 24, 32$. Only the first 10 are relevant (see Table 6). These are handled in Table 7.

In the case $A_3 \times D_5$, the non-special orbits in type $D_5$ are $1^3 2^2 3$ and $1^1 2^2 5$. This case is handled in Table 8.

| Basic Orbit | $\delta_{[w]} = d_O$ | Basic Orbit | $\delta_{[w]} = d_O$ | Basic Orbit | $\delta_{[w]} = d_O$ |
|---|---|---|---|---|---|
| $E_8$ | 0 | $E_7(a_2)$ | 8 | $(A_5 + A_1)'$ | 19 |
| $E_8(a_1)$ | 1 | $E_6 + A_1$ | 9 | $D_5(a_1) + A_2$ | 19 |
| $E_8(a_2)$ | 2 | $D_7(a_1)$ | 9 | $A_4 + A_3$ | 20 |
| $E_7 + A_1$ | 3 | $D_8(a_3)$ | 10 | $2A_3$ | 26 |
| $D_8$ | 4 | $A_7$ | 11 | $D_4 + A_1$ | 28 |
| $E_7$ | 4 | $D_6$ | 12 | $A_1 + A_2 + A_3$ | 29 |
| $E_7(a_1) + A_1$ | 5 | $D_5 + A_2$ | 13 | $2A_2 + 2A_1$ | 36 |
| $D_8(a_1)$ | 6 | $2A_4$ | 16 | $4A_1$ | 56 |
| $D_7$ | 7 | $A_5 A_2$ | 17 | | |
| $E_7(a_2) + A_1$ | 7 | $D_6(a_2)$ | 18 | | |
| $A_8$ | 8 | $A_5 + 2A_1$ | 18 | | |

Table 6. Basic orbits in $E_8$

| Orbit | $d_O$ | Character | $j\chi$ |
|---|---|---|---|
| $1^1 2^2 11$ | 4 | $((2,6),())$ | $84_x$ |
| $1^1 4^2 7$ | 7 | $((1),(3,4))$ | $400_z$ |
| $1^3 2^2 9$ | 8 | $((1,2,5),())$ | $1344_x$ |
| $3^1 4^2 5$ | 10 | $((1^2),(3^2))$ | $1050_{10}$ ⚠ |
| $1^3 4^2 5$ | 11 | $((1),(1^1 3^2))$ | $1400_{zz}$ |
| $1^1 2^4 7^1$ | 12 | $((),(2^2 4))$ | $972_x$ |
| $1^1 2^2 3^2 5$ | 13 | $((1),(2^2 3))$ | $4536_{13}$ |
| $1^3 2^4 5$ | 18 | $((),(1^1 2^2 3))$ | $4200_y$ |
| $1^3 2^2 3^3$ | 19 | $((1),(1^1 2^3))$ | $1344_w$ |
| $1^5 2^4 3$ | 26 | $((),(1^2, 2^3))$ | $840_{26} = 840_x \otimes \text{sign}$ |

Table 7. Case $D_8$

In the case $A_2 \times E_6$, the non-special orbits in type $E_6$ are $A_5$, $A_3 + A_1$, $2A_2 + A_1$, $3A_1$. This case is handled in Table 9.

There are 10 non-special nilpotent orbits in type $E_7$. The case $A_1 \times E_7$ is handled in Table 10. The left four rows represent orbits on which the $A_1$ factor has Jordan type $1^2$; the right four rows represent orbits on which the $A_1$ factor has Jordan type $2^1$.

5.6.6. *Exceptions in type $E_8$.* We now handle the exceptions marked with ⚠. For the reader's convenience, we reproduce them in Table 11

We handle the case of the orbit $3^1 4^2 5$ in type $D_8$. From Table 6, the only basic orbit with $d_O = 10$ is $D_8(a_3)$. Suppose for the sake of contradiction that

(14) $$\text{KL}^{D_8}_{E_8}(3^1 4^2 5) = \text{KL}_{E_8}(D_8(a_3)) = [w].$$

| Orbit | $d_O$ | Character | $j\chi$ |
|---|---|---|---|
| $1^3 2^2 3^1 \times 4$ | 4 | $(2^1 3^1,) \times 4$ | $84_x$ |
| $1^3 2^2 3^1 \times 1^1 3^1$ | 5 | $(2^1 3^1,) \times 1^1 3^1$ | $560_z$ |
| $1^3 2^2 3^1 \times 2^2$ | 6 | $(2^1 3^1,) \times 2^2$ | $700_x$ |
| $1^3 2^2 3^1 \times 2^1 1^2$ | 7 | $(2^1 3^1,) \times 2^1 1^2$ | $400_z$ |
| $1^1 2^2 5^1 \times 4$ | 8 | $(1^1 2^2,) \times 4$ | $1344_x$ |
| $1^1 2^2 5^1 \times 1^1 3^1$ | 9 | $(1^1 2^2,) \times 1^1 3^1$ | $3240_z$ |
| $1^3 2^2 3^1 \times 1^4$ | 10 | $(2^1 3^1,) \times 1^4$ | $1050_{10} = 1050_x$ ⚠ |
| $1^1 2^2 5^1 \times 2^2$ | 10 | $(1^1 2^2,) \times 2^2$ | $2268_x$ |
| $1^1 2^2 5^1 \times 2^1 1^2$ | 11 | $(1^1 2^2,) \times 2^1 1^2$ | $1440_{zz}$ |
| $1^1 2^2 5^1 \times 1^4$ | 14 | Not a value of $\delta_{[w]}$ | - |

TABLE 8. Case $A_3 \times D_5$

| Orbit | $d_O$ | Character | $j\chi$ |
|---|---|---|---|
| $A_5 \times 3$ | 4 | $15_q \times 3$ | $84_x$ |
| $A_3 + A_1 \times 3$ | 8 | $60_s \times 3$ | $1344_x$ |
| $2A_2 + A_1 \times 3$ | 9 | $10_s \times 3$ | $448_z$ |
| $3A_1 \times 3$ | 16 | $15_q^* \times 3$ | $3200_x$ |
| $A_5 \times 1^1 2^1$ | 5 | $15_q \times 1^1 2^1$ | $560_z$ |
| $A_3 + A_1 \times 1^1 2^1$ | 9 | $60_s \times 1^1 2^1$ | $3240_z$ |
| $2A_2 + A_1 \times 1^1 2^1$ | 10 | $10_s \times 1^1 2^1$ | $2240_x$ |
| $3A_1 \times 1^1 2^1$ | 17 | $15_q^* \times 1^1 2^1$ | $7168_w$ |
| $A_5 \times 1^3$ | 7 | $15_q \times 1^3$ | $400_z$ |
| $A_3 + A_1 \times 1^3$ | 11 | $60_s \times 1^3$ | $1400_{zz}$ |
| $2A_2 + A_1 \times 1^3$ | 12 | $10_s \times 1^3$ | $175_x$ ⚠ |
| $3A_1 \times 1^3$ | 19 | $15_q^* \times 1^3$ | $2016_w$ |

TABLE 9. Case $A_2 \times E_6$

Consider the orbit $(4^4)^+ \in \underline{\mathcal{N}_{D_8}}$. It corresponds to the character $(2^2, 2^2)^+$ of $W(D_8)$. We have $j^{W_{E_8}}_{W_{D_8}} (2^2, 2^2)^+ = 972_x$. Note that $972_x$ is the Springer representation of the nilpotent orbit $D_6 \in \underline{\mathcal{N}_{E_8}}$. By [8, Theorem 14]
$$\mathrm{KL}^{D_8}_{E_8}((4^4)^+) = \mathrm{KL}_{E_8}(D_6) = [w'].$$
Moreover, $54^2 3 > (4^4)^+$ in the closure ordering.

**Definition 49** ([24]). *Let $[w_1]$, $[w_2]$ be conjugacy classes in $W$. We say that $[w_1] > [w_2]$ if $\overline{(L^\heartsuit \mathfrak{c})_{[w_1]}} \supset (L^\heartsuit \mathfrak{c})_{[w_2]}$.*

By Proposition 18, we have $[w] > [w']$. On the other hand, by Proposition 18 again, this implies that $D_6 < D_8(a_3)$. But this is known to be false: see the closure graph [23, Page 247]. Note that there $D_8(a_3)$ is denoted $E_8(b_6)$. Therefore (14) does not hold.

| Orbit $\times 1^2$ | $d_O$ | Character $\times 1^2$ | $j\chi$ | Orbit $\times 2^1$ | $d_O$ | Character $\times 2^1$ | $j\chi$ |
|---|---|---|---|---|---|---|---|
| $D_6$ | 5 | $35_4 = 35_b$ | $560_z$ | $D_6$ | 4 | $35_b$ | $84_x$ |
| $D_6(a_2)$ | 9 | $280_8 = 280_b$ | $3240_z$ | $D_6(a_2)$ | 8 | $280_8$ | $1344_x$ |
| $(A_5)'$ | 10 | $216_9 = 216'_a$ | $1050_x$ ⚠ | $(A_5)'$ | 9 | $216_9$ | $3240_z$ |
| $(A_5 + A_1)''$ | 10 | $70_9 = 70'_a$ | $2240_x$ | $(A_5 + A_1)''$ | 9 | $70_9$ | $448_z$ |
| $D_4 + A_1$ | 13 | $84_{12} = 84_a$ | $4536_z$ | $D_4 + A_1$ | 12 | $84_{12}$ | $972_x$ |
| $A_3 + 2A_1$ | 17 | $216_{16} = 216_a$ | $7168_w$ | $A_3 + 2A_1$ | 16 | $216_{16}$ | $3200_x$ |
| $(A_3 + A_1)'$ | 18 | $280_{17} = 280_b^*$ | $4200_y$ | $(A_3 + A_1)'$ | 17 | $280_{17}$ | $7168_w$ |
| $2A_2 + A_1$ | 19 | $70_{18} = 70_a$ | $2016_w$ | $2A_2 + A_1$ | 18 | $70_{18}$ | $3150_y$ |
| $4A_1$ | 29 | $15_{28} = 15_a$ | $1400_{29} = 1400_{zz}^*$ | $4A_1$ | 28 | $15_{28}$ | $700_{28}$ |
| $(3A_1)'$ | 32 | Not a value of $\delta_{[w]}$ | - | $(3A_1)'$ | 31 | Not a value of $\delta_{[w]}$ | - |

TABLE 10. Case $A_1 \times E_7$

| Case | Orbit | $d_O$ | Character | $j\chi$ |
|---|---|---|---|---|
| $D_8$ | $3^1 4^2 5$ | 10 | $(1^2, 3^2)$ | $1050_{10}$ |
| $A_3 \times D_5$ | $1^1 2^2 5^1 \times 1^4$ | 10 | $(2^1 3^1,) \times 1^4$ | $1050_{10}$ |
| $E_7 \times A_1$ | $(A_5)' \times 1^2$ | 10 | $216'_a \times 1^2$ | $1050_x$ |
| $A_2 \times E_6$ | $2A_2 + A_1 \times 1^3$ | 12 | $10_s \times 1^3$ | $175_x$ |

TABLE 11. Orbits such that $j\chi$ is not Springer

We handle the case of the orbit $O' = 1^1 2^2 5^1 \times 1^4$ in type $A_3 \times D_5$. Again, $d_{O'} = 10$, so suppose for the sake of contradiction that

(15) $\qquad \mathrm{KL}_{E_8}^{A_3 \times D_5}(1^1 2^2 5^1 \times 1^4) = \mathrm{KL}_{E_8}(D_8(a_3)) = [w]$.

Consider the nilpotent orbit $5^1 1^5$ in type $D_5$. It symbol is $(, 3^1 1^2)$. From Alvis' tables, $j((, 3^1 1^2), 1^4) = 972_x$ which corresopnds to the orbit $D_6$. Moreover, $5^1 1^5 \times 1^4 < 5^1 2^2 1^1 \times 1^4$. Using the same reasoning as before, (15) cannot hold.

We handle the case of the orbit $O' = (A_5)' \times 1^2$ in type $E_7 \times A_1$. Again, $d_{O'} = 10$, so suppose for the sake of contradiction that

(16) $\qquad \mathrm{KL}_{E_8}^{E_7 \times A_1}((A_5)' \times 1^2) = \mathrm{KL}_{E_8}(D_8(a_3)) = [w]$.

Consider the orbit $O'' = (A_5)'' \times 1^2$ in type $E_7 \times A_1$. This orbit satisfies $d_{O''} = 13$, and $O'' \not\subset \overline{O'}$. Moreover, $O''$ is special with symbol $105_c \times 1^2$. Its $j$-induction to $W(E_8)$ is $j_{W(E_7 \times A_1)}^{W(E_8)} 105_c \times 1^2 = 2800_z$ which is the Springer representation associated to the orbit $D_5 + A_2$. Let

$$[w'] = \mathrm{KL}_{E_8}^{E_7 \times A_1}((A_5)'' \times 1^2) = \mathrm{KL}_{E_8}(D_5 + A_2)$$

where the last equality is by [8, Theorem 14]. Since $D_5 + A_2 \subset \overline{D_8(a_3)}$, Proposition 18 says that $[w'] < [w]$. On the other hand, if (16) held, then Proposition 18 would show that $O'' \subset \overline{O'}$, which is a contradiction.

Finally we handle the case of the orbit $O' = 2A_2 + A_1 \times 1^3$ in type $E_6 \times A_2$. From Table 6, there is only one elliptic $[w]$ with $\delta_{[w]} = d_{O'} = 12$, so suppose for the sake of contradiction that

$$\text{KL}^{E_6 \times A_2}_{E_8}(2A_2 + A_1 \times 1^3) = \text{KL}_{E_8}(D_6) = [w]. \tag{17}$$

Consider the orbit $O'' = A_3 \times 1^3$ in type $E_6 \times A_2$. This orbit satisfies $d_{O''} = 13$, and $O'' \not\subset \overline{O'}$. Moreover, $O''$ is special with symbol $81_{10} \times 1^3$. Its $j$-induction to $W(E_8)$ is $j^{W(E_8)}_{W(E_6 \times A_2)} 81_{10} \times 1^3 = 4536_z$, which is the Springer representation associated to the orbit $D_5 + A_2$. Let

$$[w'] = \text{KL}^{E_6 \times A_2}_{E_8}(A_3 \times 1^3) = \text{KL}_{E_8}(D_5 + A_2)$$

where the last equality is by [8, Theorem 14]. Since $D_5 + A_2 \subset \overline{D_6}$, Proposition 18 says that $[w'] < [w]$. On the other hand, if (17) held, then Proposition 18 would show that $O'' \subset \overline{O'}$, which is a contradiction.

This finishes the proof in type $E_8$.